\documentclass[a4paper]{article}
\usepackage{amsmath,amsthm,amsfonts,thmtools}
\usepackage{amssymb}
\usepackage{newtxtext}
\usepackage[cmintegrals]{newtxmath}
\usepackage{graphicx}
\usepackage{subfigure}
\usepackage{euscript}
\usepackage{units}
\usepackage{enumerate}
\usepackage{enumitem}
\usepackage{xcolor}
\usepackage[bookmarks=true,hidelinks]{hyperref}
\usepackage{bbm}

\numberwithin{equation}{section}
\newtheorem{theorem}{Theorem}[section]

\newtheorem{lemma}[theorem]{Lemma}

\newtheorem{algorithm}[theorem]{Algorithm}

\numberwithin{equation}{section}

\theoremstyle{definition}

\newtheoremstyle{myremarkstyle}{}{}{}{}{\bfseries}{.}{ }{}
\theoremstyle{myremarkstyle}
\declaretheorem[name=Remark,qed=$\blacksquare$,numberlike=theorem]{remark}

\makeatletter
\newcommand*{\intavg}{%
  \mint@l{-}{}%
}
\newcommand*{\mint@l}[2]{%
  \@ifnextchar\limits{%
    \mint@l{#1}%
  }{%
    \@ifnextchar\nolimits{%
      \mint@l{#1}%
    }{%
      \@ifnextchar\displaylimits{%
        \mint@l{#1}%
      }{%
        \mint@s{#2}{#1}%
      }%
    }%
  }%
}
\newcommand*{\mint@s}[2]{%
  \@ifnextchar_{%
    \mint@sub{#1}{#2}%
  }{%
    \@ifnextchar^{%
      \mint@sup{#1}{#2}%
    }{%
      \mint@{#1}{#2}{}{}%
    }%
  }%
}
\def\mint@sub#1#2_#3{%
  \@ifnextchar^{%
    \mint@sub@sup{#1}{#2}{#3}%
  }{%
    \mint@{#1}{#2}{#3}{}%
  }%
}
\def\mint@sup#1#2^#3{%
  \@ifnextchar_{%
    \mint@sub@sup{#1}{#2}{#3}%
  }{%
    \mint@{#1}{#2}{}{#3}%
  }%
}
\def\mint@sub@sup#1#2#3^#4{%
  \mint@{#1}{#2}{#3}{#4}%
}
\def\mint@sup@sub#1#2#3_#4{%
  \mint@{#1}{#2}{#4}{#3}%
}
\newcommand*{\mint@}[4]{%
  \mathop{}%
  \mkern-\thinmuskip
  \mathchoice{%
    \mint@@{#1}{#2}{#3}{#4}%
        \displaystyle\textstyle\scriptstyle
  }{%
    \mint@@{#1}{#2}{#3}{#4}%
        \textstyle\scriptstyle\scriptstyle
  }{%
    \mint@@{#1}{#2}{#3}{#4}%
        \scriptstyle\scriptscriptstyle\scriptscriptstyle
  }{%
    \mint@@{#1}{#2}{#3}{#4}%
        \scriptscriptstyle\scriptscriptstyle\scriptscriptstyle
  }%
  \mkern-\thinmuskip
  \int#1%
  \ifx\\#3\\\else_{#3}\fi
  \ifx\\#4\\\else^{#4}\fi  
}
\newcommand*{\mint@@}[7]{%
  \begingroup
    \sbox0{$#5\int\m@th$}%
    \sbox2{$#5\int_{}\m@th$}%
    \dimen2=\wd0 %
    \let\mint@limits=#1\relax
    \ifx\mint@limits\relax
      \sbox4{$#5\int_{\kern1sp}^{\kern1sp}\m@th$}%
      \ifdim\wd4>\wd2 %
        \let\mint@limits=\nolimits
      \else
        \let\mint@limits=\limits
      \fi
    \fi
    \ifx\mint@limits\displaylimits
      \ifx#5\displaystyle
        \let\mint@limits=\limits
      \fi
    \fi
    \ifx\mint@limits\limits
      \sbox0{$#7#3\m@th$}%
      \sbox2{$#7#4\m@th$}%
      \ifdim\wd0>\dimen2 %
        \dimen2=\wd0 %
      \fi
      \ifdim\wd2>\dimen2 %
        \dimen2=\wd2 %
      \fi
    \fi
    \rlap{%
      $#5%
        \vcenter{%
          \hbox to\dimen2{%
            \hss
            $#6{#2}\m@th$%
            \hss
          }%
        }%
      $%
    }%
  \endgroup
}

\renewcommand{\geq}{\geqslant}
\renewcommand{\leq}{\leqslant}

\renewcommand{\epsilon}{\varepsilon}
\renewcommand{\phi}{\varphi}

\newcommand{\R}{\mathbb{R}}

\newcommand{\Dx}{{\Delta x}}

\newcommand{\Dt}{\Delta t}

\newcommand{\hf}{{\unitfrac{1}{2}}}

\newcommand{\jphf}{{j+\hf}}
\newcommand{\jmhf}{{j-\hf}}

\begin{document}

\date{\today}

\title{A machine learning framework for data driven acceleration of computations of differential equations.}

\author{Siddhartha Mishra \thanks{Seminar for Applied Mathematics (SAM), D-Math \newline
  ETH Z\"urich, R\"amistrasse 101, 
  Z\"urich-8092, Switzerland}}

\date{\today}

\maketitle
\begin{abstract}
We propose a machine learning framework to accelerate numerical computations of time-dependent ODEs and PDEs. Our method is based on recasting (generalizations of) existing numerical methods as artificial neural networks, with a set of trainable parameters. These parameters are determined in an offline training process by (approximately) minimizing suitable (possibly non-convex)
loss functions by (stochastic) gradient descent methods. The proposed algorithm is designed to be always consistent with the underlying differential equation. Numerical experiments involving both linear and non-linear ODE and PDE model
problems demonstrate a significant gain in computational efficiency over standard numerical methods.
\end{abstract}

\section{Introduction}
Differential equations, both ordinary and partial, are ubiquitous in science and engineering. It is not possible to obtain explicit solution formulas for differential equations, except in the simplest cases.
Hence, numerical approximations of differential equations constitutes a key tool in their study. A wide variety of numerical methods have been developed to approximate differential equations robustly and
efficiently. For the initial value problem for ordinary differential equations, popular methods include Runge-Kutta and multi-step methods, \cite{HW1,LEV2} and references therein. Widely used numerical methods for approximating
PDEs include finite difference \cite{LEV2}, finite volume \cite{GR1}, finite element \cite{BS1} and spectral \cite{TREF1} methods. 

The exponential increase in computational power in the last decades provides the opportunity for solving very challenging, large scale computational problems for differential equations, such as uncertainty
quantification (UQ) \cite{UQbook,UQhb}, (Bayesian) inverse problems \cite{STU1} and (real time) optimal control, design and constrained optimization \cite{OC,OC1}. One requires
very large number of fast approximations of ODE and PDEs to solve such problems, for instance when evaluating Monte Carlo samples in an UQ or Bayesian inverse problem framework. Currently available numerical methods, particularly for nonlinear PDEs, tend to be too slow to allow such realistic computations.

Machine learning, in the form of artificial neural networks (ANNs), has become extremely popular in computer science in recent years. This term is applied to a plethora of methods that aim to approximate functions with layers of
units (neurons), connected by linear operations between units and nonlinear activations within units, \cite{DLbook} and references therein. \emph{Deep learning}, i.e an artificial neural network with a large number of intermediate (hidden) layers has proven extremely successful at diverse tasks, for instance in image processing, signal processing and natural language processing \cite{DL-nat}. A key element in deep learning is the \emph{training} of tuning parameters in the underlying neural network by (approximately) minimizing suitable \emph{loss functions}. The resulting (non-convex) optimization problem, on a very high dimensional parameter space, can be efficiently solved with variants of the stochastic gradient descent method \cite{ADAM,SG}.

Machine learning methods, particularly deep learning, are being increasingly used in the context of numerical computation of differential equations. As this is a rapidly evolving field, we will only attempt a skeletal literature survey here. One class of methods attempt to replace numerical schemes for differential equations by deep networks, see \cite{WE1,WE2,KAR1,JR1} and references therein. These methods have been successfully used in different contexts, for instance in approximating very high dimensional problems arising in mathematical finance \cite{WE1}, by exploiting integral representation formulas for the underlying solutions. However, it is as yet unclear if an end to end deep neural network can learn the physics of the underlying PDE in the absence of such formulas. This could constitute a stumbling block in approximating solutions of complicated nonlinear PDEs with deep learning. 

Another school of thought aims to augment existing numerical methods by embedding deep learning modules within them. As examples, one can think of solving the pressure Poisson equation within an incompressible flow solver by a convolutional neural network as in \cite{INC} or learning troubled cell indicators in a RKDG code by a deep network as in \cite{DR1}. 

In this paper, we propose a variant of the machine learning framework for approximating (time-dependent) differential equations. Our starting point is the observation that evaluating approximate solutions of ODEs and PDEs on coarse space-time grids is very cheap computationally.  However, the accuracy of such coarse grid representations is rather poor. Consequently, we will use a machine learning framework to train explicit or implicit parameters in (generalizations of) standard numerical methods in order to minimize a loss (error) function that measures the difference of the (trained) solution on coarse grids with projections of fine grid solutions. The resulting scheme will hopefully be significantly more accurate than the underlying standard method on the coarse grid, while being as computationally cheap. We motivate our general strategy with the following simple example,
\begin{figure}[htbp]
\centering
\includegraphics[width=0.45\linewidth]{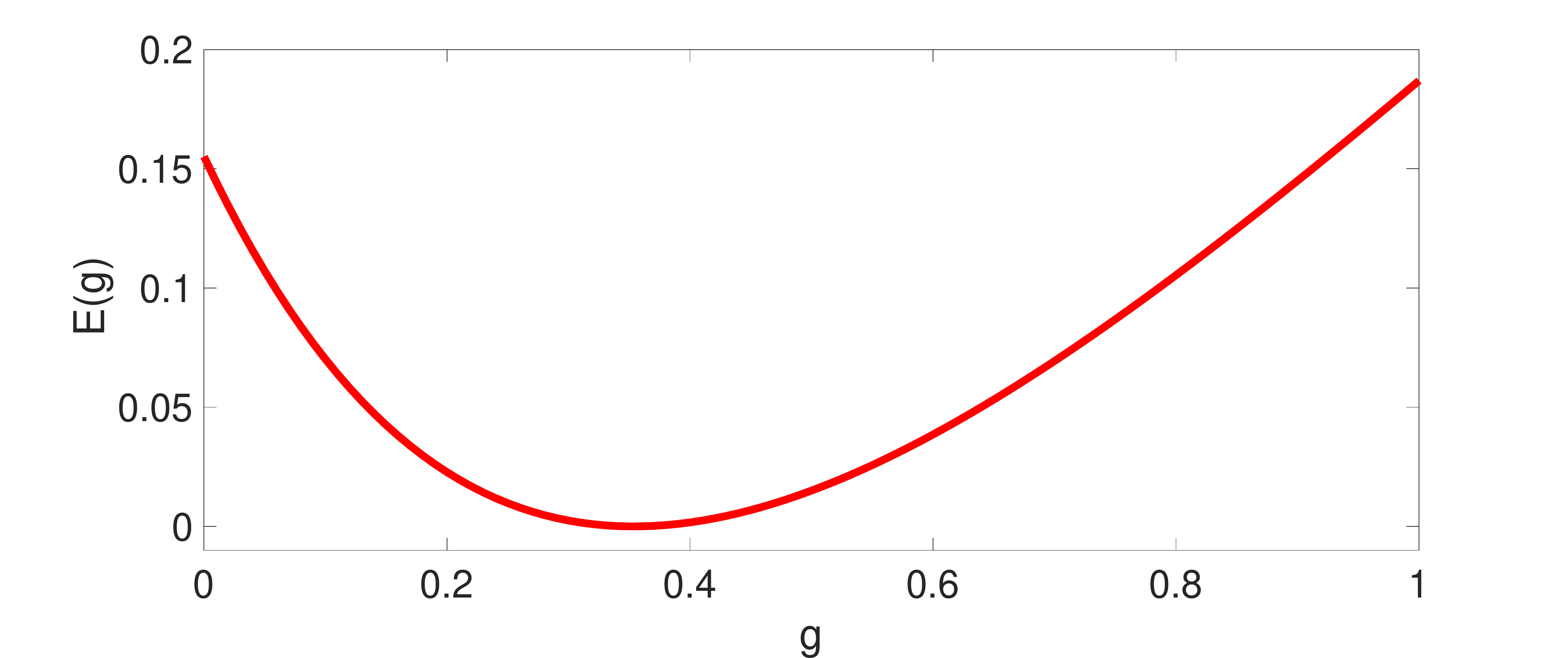} \includegraphics[width=0.45\linewidth]{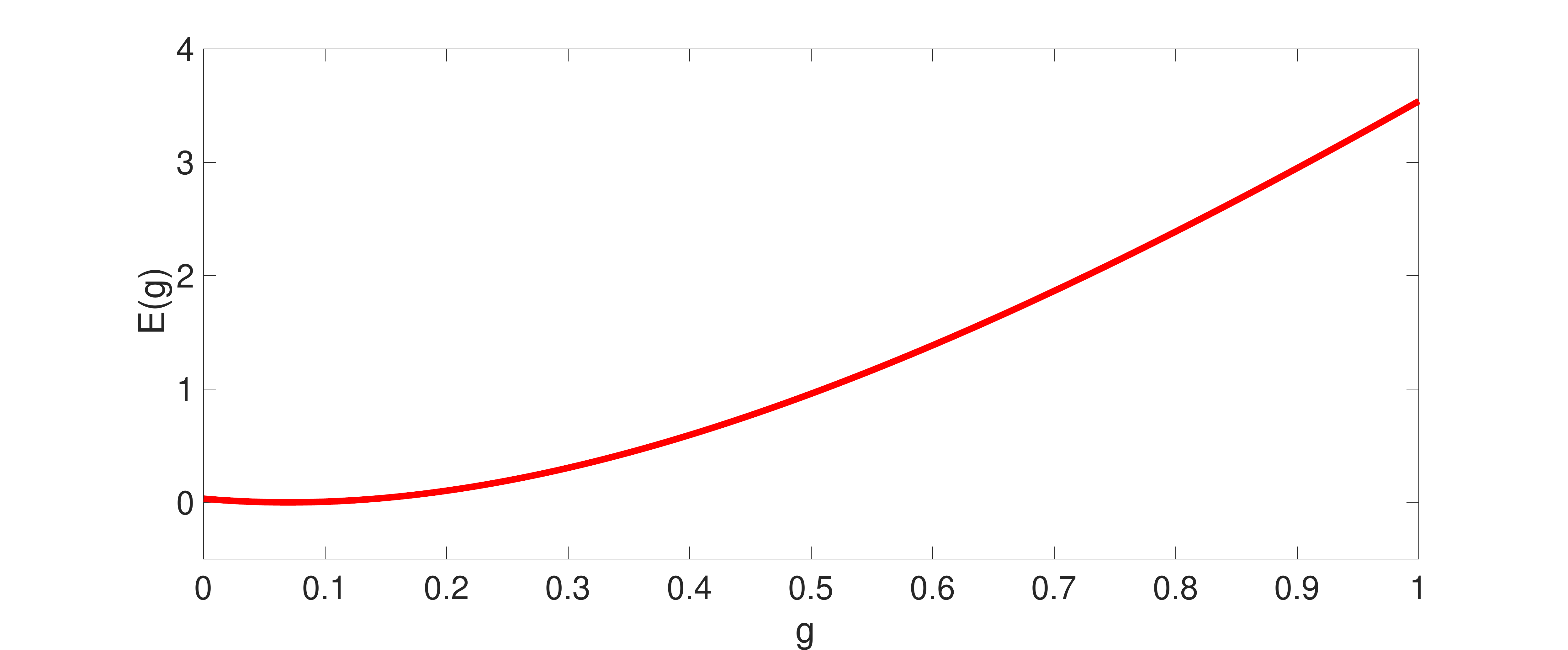}
\caption{Error $E_2(g)$ \eqref{eq:bdfge} (Y-axis) vs. the scheme parameter $g$ (X-axis) in the generalized BDF scheme \eqref{eq:bdfg1} for $\Dt = 0.5$ and two different values of $c$. Left $c=1$, the minimum in the error is achieved at
$g = 0.35$ and yields a factor of $39.97$ times reduction in the error over the standard three-point BDF2 scheme i.e \eqref{eq:bdfg} with $g = 0.5$. Right: $c=5$: The minimum is at $g=0.07$ and the error is reduced by a factor of $677.95$ over the standard BDF2 scheme}
\label{fig:1}
\end{figure}

\subsection{A motivating example}
We consider the following autonomous ODE,
\begin{equation}
\label{eq:ode1}
\begin{aligned}
u^{\prime}(t) = F(u(t)), \quad t \in (0,T), \\
u(0) = u_0.
\end{aligned}
\end{equation}
Here, $u:[0,T] \rightarrow \R^d$ is the vector of unknowns and $F \in {\rm Lip}_{{\rm loc}}(\R^d)$ is the right hand side. 

A very popular class of numerical methods to solve \eqref{eq:ode1}, particularly for \emph{stiff} right hand sides, are the so-called implicit multi-step methods or backward difference formulas (BDFs) \cite{HW1}. Assuming a uniform time step $\Dt$, denoting the time level by $t_n = n \Dt$ and the approximate solution by $U_n \approx u(t_n)$, a general form of a three-point BDF is given by,
\begin{equation}
\label{eq:bdfg}
(1+g_{n+2})U_{n+2} - (1+2g_{n+2})U_{n+1} + g_{n+2}U_n = \Dt F(U_{n+2}),
\end{equation}
for $n \geq 0$ and $g_{n+2} \in \R$ for each $n$. We need initial values $U_0,U_1$ to march in time in \eqref{eq:bdfg}. It is straightforward to check using Taylor expansions that for any $g_{n+2} \in \R$, the method \eqref{eq:bdfg} is consistent with \eqref{eq:ode1} and at least first-order accurate. By setting $g_{n+2}=0,~\forall n$, we recover the standard backward Euler method while $g_{n+2}=\frac{1}{2}$ for all $n$ yields the second-order accurate BDF2 method. It is customary to assume that a second-order accurate method is preferable. Hence, one always sets $g_{n+2} \equiv 0.5$. 

However, this point of view does not take the data of the specific problem that we are solving into account, namely the non-linearity $F$, the dimension $d$, the initial data $u_0$ and the grid size $\Dt$. It is not a priori obvious if the choice of $g_{n+2} = 0.5$ will provide the best solution (the least error) for a given data set. In fact, this is not true in general. As an example, we consider the simplest linear scalar ODE by setting $d=1$ and $F(u) = - cu$, for some $c \in \R_+$, in \eqref{eq:ode1}. In this case, the explicit solution is given by
\begin{equation}
\label{eq:ode1s}
u(t) = u_0 e^{-ct}.
\end{equation}
Now setting $U_0 = u_0$, $U_1 = e^{-c\Dt}u_0$, the solution computed at the second time level $2\Dt$, by the generalized form of the three point BDF scheme \eqref{eq:bdfg} is 
\begin{equation}
\label{eq:bdfg1}
U_2 = \frac{1+2g}{1+g+c\Dt} U_1 - \frac{g}{1+g+c\Dt}U_0.
\end{equation}
Hence, the local error at the second time level is
\begin{equation}
\label{eq:bdfge}
E_2(g) =|U_2 - e^{-2c\Dt}u_0|^2,
\end{equation}
It is straightforward to observe that the minimizer,
$$
g^{\ast} = {\rm arg}~\min\limits_{g \in \R} E_2(g),
$$
depends explicitly on the parameters $\Dt$ and $c$ and is not universally $g^{\ast} = 0.5$. In figure \ref{fig:1}, we plot the function $E_2(g)$ for $\Dt = 0.5$ and two different values of $c$, namely $c=1$ and $c=5$, respectively. We see from this figure that the loss (error) function is convex in the parameter $g$ and the error is vastly reduced in a minimization process, i.e, by a factor of $39.97$ for $c=1$ and a factor of $677.95$ for $c=5$, respectively. Hence, by focusing on a specific data set, we can potentially obtain speed ups of two to three orders of magnitude for this simple ODE.
\subsection{Aims and scope of this paper}
Our objective is to generalize the strategy presented in the above motivating example. We will recast (generalizations of) standard numerical methods for time-dependent ODEs and PDEs as \emph{multi-layer} artificial neural networks with a set of trainable parameters. The resulting network will always be designed to be consistent with the underlying differential equation by constraining the parameter set. During an offline training phase, we will train these parameters by (approximately) minimizing a loss function, over the parameter states, by a suitable (stochastic) gradient descent method. The efficiency of the resulting trained scheme is verified on a \emph{test set}. Thus, our methods will be (rather restricted) types of (deep) neural networks for approximating time-dependent differential equations. 

We organize the rest of the paper as follows: in section \ref{sec:2}, we present our abstract machine learning framework. In section \ref{sec:3}, we apply the proposed algorithm to ordinary differential equations. The machine learning algorithm is 
applied to the heat equation, linear transport equation, scalar conservation laws and the Euler equations of gas dynamics in sections \ref{sec:5}, \ref{sec:4}, \ref{sec:6} and \ref{sec:7}. We summarize the contents of the paper and provide a perspective in section \ref{sec:8}.

\section{The abstract machine learning framework}
\label{sec:2}
For definiteness, we consider a one-dimensional nonlinear time-dependent PDE of the form,
\begin{equation}
\label{eq:pde1}
\begin{aligned}
u_t &= L(u,u_x,u_{xx}), \quad (x,t) \in [X_l,X_r] \times [0,T], \\
u(0,x) &= u_0(x,\omega), \\
L_b u(X_l,t) &= u_l(t,\gamma_l), \quad L_bu(X_r,t) = u_r(t,\gamma_r).
\end{aligned}
\end{equation}
Here, $u:  [X_l,X_r] \times [0,T] \rightarrow \R^m$ is the vector of unknowns. $L$ is a possibly non-linear differential operator involving both the first and second spatial derivatives of $u$ (interpreted as vectors) that will be specified in
subsequent examples but is kept deliberately ambiguous here for the sake of generality. $L_b$ refers to a \emph{boundary} operator and the initial and boundary data depend on parameters $\omega,\gamma_{l,r} \in \R^D$ for possibly $D >> 1$.

For simplicity, we discretize $[X_l,X_r]$ on a uniform grid with mesh size $\Dx$ and denote the discrete points as $x_j = X_l + j Dx$, for $0 \leq j \leq J+1$. Similarly, we choose a uniform time step $\Dt$ and denote the $n$-th time level as $t^n = n\Dt$, with $T = N\Dt$, and denote the approximate solution (by a finite difference scheme as ) as $U^n_j \approx u(x_j,t^n)$. Moreover, one can readily consider a finite volume or finite element discretization by letting $U^n_j$ approximate cell averages or nodal point values, respectively.  We denote the vector, $U^n = \{ U^n_j \}_{1\leq j \leq N}$.

On this grid, we discretize the abstract PDE \eqref{eq:pde1} with the following numerical method, 
\begin{equation}
\label{eq:abnm}
U^{n+1} = L^n(\Dt,\Dx, U^n, A^n U^n, F(U_n)), R(U^n)).
\end{equation}

Several remarks are in order about the form of the abstract scheme \eqref{eq:abnm}. First, $A^n$ is a linear operator that subsumes potentially multiple applications of matrices. Second, the function $F$ denotes a
composite of nonlinearities that occur within the non-linear operator $L$ in \eqref{eq:pde1}. Furthermore, we write \eqref{eq:abnm} as a one-step explicit scheme. However, we implicitly assume that even if the underlying time discretization was (semi-)implicit, the resulting system of nonlinear equations is solved by an iterative procedure, such as a Newton method, and the results of the Newton steps are subsumed in the general form of the right hand side term $L^n$ in \eqref{eq:abnm}. In particular, such iterations might involve multiple applications of the non-linearity $F$ in \eqref{eq:abnm} and could involve the function $R$ in \eqref{eq:abnm} that might include (multiple compositions of) the standard ReLU function,
\begin{equation}
\label{eq:relu}
\sigma (w) = \max(w,0).
\end{equation}

\begin{figure}[htbp]
\centering
\includegraphics[width=8cm]{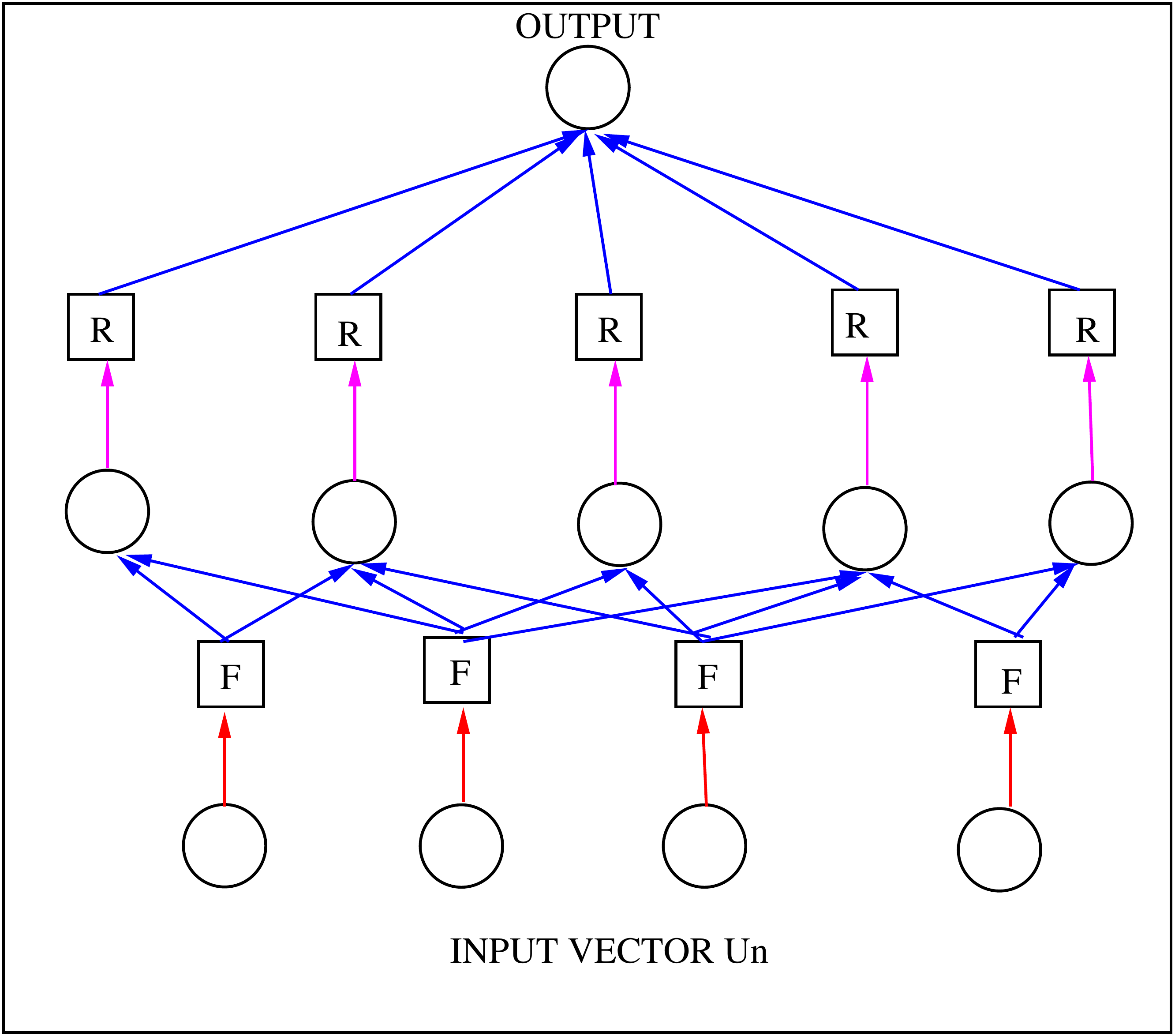}
\caption{A simplified representation of a single time step of the abstract numerical scheme \eqref{eq:abnm} as a \emph{multi-layer neural network}. We show a part (4 components) of the solution vector $U^n$ in \eqref{eq:abnm} as the input. The blue arrows represent linear mappings ($A^n$) between units. The magenta arrows represent nonlinear (scalar) activations, for instance by the standard ReLU function \eqref{eq:relu} of units and the red arrows represent evaluation of nonlinearities $F$, that are inherent to the underlying differential equation \eqref{eq:pde1}. The hidden layers might correspond to steps of a Newton method. The output is a single component of the vector $U^{n+1}$ in \eqref{eq:abnm}. Note that the function $F$ can be replaced by an artificial neural network module to realize the whole scheme \eqref{eq:abnm} as an artificial neural network in the sense of \cite{DLbook}.}
\label{fig:2}
\end{figure}

The whole method can be represented as a \emph{neural network} as shown in figure \ref{fig:2}, but with non-linearities $F$, based on the underlying PDE instead of scalar activation functions as in a traditional deep network \cite{DLbook}. However, given the universal approximation property of standard artificial neural networks \cite{BAR1,HOR1}, one can approximate the underlying nonlinearities $F$ in \eqref{eq:abnm} by artificial neural networks. In particular, networks with a few hidden layers can be trained to approximate smooth nonlinear functions \cite{ YAR1}. Hence, one can think of the module $F$ in \eqref{eq:abnm1}, figure \ref{fig:2} as an \emph{additional neural network}, realizing the whole scheme \eqref{eq:abnm} as a neural network in the sense of \cite{DLbook}. However, for the sake of consistency and computational efficiency, we perform a direct evaluation of the nonlinearity $F$ in \eqref{eq:abnm} in this paper. A very concrete realization of the neural network representation for numerical schemes is provided in section \ref{sec:6} for a Rusanov type scheme approximating scalar conservation laws, see figure \ref{fig:7}.

We constrain the numerical method (or alternatively the artificial neural network) \eqref{eq:abnm}  to be consistent with the PDE \eqref{eq:pde1} by imposing constraints on the linear operator $A^n$ (and the structure of the neural network approximating $F$). Further stability conditions on the scheme can also be imposed. Consequently, we can rewrite the numerical method \eqref{eq:abnm} in the following parametric form,
\begin{equation}
\label{eq:abnm1}
U^{n+1} = L^n(\Dt,\Dx,U^n,\theta^n),
\end{equation}
 in terms of a parameter vector $\theta^n \in \R^d$. We impose the constraints on the scheme such that by choosing $\theta^n = \overline{\theta}^n$, leads to the recovery of standard numerical methods, for instance the choice of $g = 0.5$ in the scheme \eqref{eq:bdfg} yields the standard second-order three-point BDF2 scheme for discretizing the ODE \eqref{eq:ode1}.

For generating the training set, we select a certain subset of the parameter space $\{\omega_i,\gamma_{1,i}, \gamma_{2,i}\}_{1\leq i \leq M}$ where $M >> 1$ is the size of the training set. Let $\Dt_f << \Dt$ and $\Dx_f << \Dx$ be the time step and mesh size for a (very) fine grid (uniform) discretization of $[0,T] \times [X_l,X_r]$. For training data 
$\{ u_0(\omega_i), u_l(\gamma_{1,i}), u_r(\gamma_{2,i})\}_{i}$, we approximate \eqref{eq:pde1} with the scheme \eqref{eq:abnm1}, with a particular choice of $\theta^n = \overline{\theta}^n$, on the fine grid.  The resulting solution, denoted as $U_{{\rm ref}}^{n,i}$, is obtained by projecting the fine grid solution to the coarse grid, either with cell averaging or point wise sampling. This training data is generated offline and will be rather computationally expensive as a fine grid solution needs to be computed. 

Next we set up a \emph{training loss function} by defining the error,
\begin{equation}
\label{eq:lf1}
E(\theta) := \frac{1}{p} \sum\limits_{i}\sum\limits_{n=1}^{N} \|U^{n,i} (\theta^n)- U^{n,i}_{{\rm ref}} \|_{L^p}^p.
\end{equation}
Here, $1 \leq p < \infty$ and we usually consider $p=1$ or $p=2$ and denote $\theta = \{ \theta^n \}_{1 \leq n \leq N}$ as the combined vector of trainable parameters. 

The objective of the training process is to minimize the loss function \eqref{eq:lf1} i.e, find a minimizer $\theta^{\ast}$:
\begin{equation}
\label{eq:lf2}
\theta^{\ast} = {\rm arg}~\min\limits_{\theta \in \R^{Nd}} E(\theta).
\end{equation}
The resulting \emph{trained} time-marching scheme is
\begin{equation}
\label{eq:abnm2}
U^{n+1} = L^n(\Dt,\Dx,U^n,\theta^{n,\ast}).
\end{equation}

We summarize the resulting algorithm for our machine learning framework below:
\begin{algorithm}
\label{alg:1}
Given a specific model i.e, underlying ODE or PDE, for instance \eqref{eq:pde1} on a specific space-time domain, we 
\begin{description}
\item [{\textbf Step $1$}:] Choose a consistent (and stable) numerical method (alternatively neural network) for approximating the underlying differential equation, for instance the general form \eqref{eq:abnm} or \eqref{eq:abnm1}. This numerical method will approximate solutions of \eqref{eq:pde1} on a coarse grid.
\item [{\textbf Step $2$}:] Generate the \emph{training set} by choosing a specific (finite but possibly large) parametric set of initial (and boundary conditions) and approximating the underlying PDE with a standard numerical method, for
instance \eqref{eq:abnm1} but with the parameter set $\overline{\theta}$, on a (very) fine grid. Then, project the fine grid solutions on the coarse grid to create the training set. 
\item [{\textbf Step $3$}:]  Set up the loss function \eqref{eq:lf1} and use a gradient descent method to (approximately) find a local minimum of this possibly non-convex function over the parameter space. The gradient descent method can be initialized with the parameter values $\overline{\theta}$, corresponding to a standard numerical method. 
\item [{\textbf Step $4$}:]  The minimizers $\theta^{\ast}$ serve as the parameters in the trained scheme \eqref{eq:abnm2}. The trained scheme is run on a \emph{test set} in the online phase, to ascertain gains in computational efficiency.

\end{description}
\end{algorithm}

We remark that the algorithm \ref{alg:1} is guaranteed to reduce the error, on the underlying coarse grid, over a standard scheme (i.e \eqref{eq:abnm1} with parameter set $\overline{\theta}$), on the training set. Moreover, the trained scheme \eqref{eq:abnm2} will always be a consistent (and stable) discretization of the underlying model. It is difficult to obtain theoretical guarantees on the amplitude of the gain in computational efficiency for our machine learning framework, on the \emph{test set}. This gain depends on the underlying model, numerical method, grid size, choice of training set and on the efficacy of the gradient descent in finding a minimum. We will test this machine learning algorithm on a variety of problems in the following sections in order to empirically demonstrate its efficiency. 

\section{ODEs}
\label{sec:3}
In this section, we will test algorithm \ref{alg:1} on the following two ordinary differential equations,
\subsection{A linear ODE}
We start with the following second-order linear ODE modeling oscillators, 
\begin{equation}
\label{eq:odel}
u^{\prime \prime}(t) + c^2 u(t) = 0, \quad u(0) = u_0.
\end{equation}
We can readily write \eqref{eq:odel} as a first-order system by introducing the auxiliary variable $v = u^{\prime}$, resulting in 
\begin{equation}
\label{eq:odel1}
\begin{aligned}
u^{\prime} &= - cv, \quad 
v^{\prime} &= cu.
\end{aligned}
\end{equation}
It is easy to see that \eqref{eq:odel} (equivalently \eqref{eq:odel1}) has  an explicit solution given by,
\begin{equation}
\label{eq:odel2}
u(t) = u_0 \cos(ct), \quad v(t) = u_0\sin(ct).
\end{equation}

For the sake of simplicity, we choose the
generalized three-point backward difference formula \eqref{eq:bdfg}, with $U = [u,v]$ and $F(U) = [-cv, cu]$, as the underlying numerical scheme i.e step $1$ of Algorithm \ref{alg:1}. We choose a uniform grid in time with time step $\Dt$ and initialize the scheme \eqref{eq:bdfg}
with $U_0 = [u_0,u_0]$ and $U_1 = [u_0\cos(c\Dt),u_0\sin(c\Dt)]$.

Our objective is to approximate the solution of \eqref{eq:odel1} at the next two time levels, i.e determine $U_2$ and $U_3$ by \eqref{eq:bdfg} with $\Dt = \frac{1}{3}$.  As the constant $c$ determines the frequency of oscillations in \eqref{eq:odel}, a grid with a time step of $\Dt = \frac{1}{3}$ is extremely coarse for large values of $c$, as several oscillations occur within a single time step. The task at hand is to determine whether the machine learning algorithm \ref{alg:1} can (significantly) improve the accuracy of the scheme \eqref{eq:bdfg} on such a coarse grid. 

For step $2$ of Algorithm \ref{alg:1}, we generate the training set by randomly selecting $I$ data points on the interval $[0,1]$ and denoting them as $\{u_0^i\}$ with $1 \leq i \leq I=10$. Corresponding to these data points, we use the exact solutions \eqref{eq:odel2} at times $t_2= \frac{2}{3}$ and $t_3 = 1$ to generate the reference training data $U^{n,i}_{{\rm ref}}$ for all $i$ and $n=2,3$. 

In step $3$ of Algorithm \ref{alg:1}, we set up the $l^2$ error,
\begin{equation}
\label{eq:odel4}
E_2(g_2,g_3) : = \frac{1}{2} \sum\limits_{i=1}^I\sum\limits_{n=2,3} |U^{n,i}(g_n) - U^{l,i}_{{\rm ref}} |^2,
\end{equation}
and minimize the loss function \eqref{eq:odel4} over two parameters $g_{2,3} \in \R \times \R$. Given this simple two-dimensional (in parameters) problem, the minimization is performed by a standard steepest gradient descent initialized at the point $(g_2,g_3) = (0.5,0.5)$, corresponding to the \emph{second-order} BDF2 scheme. 
\begin{table}[htbp]
\centering
\begin{tabular}{|c|c|c|c|}
\hline
$c$ & $g^{\ast}_2$ & $g^{\ast}_3$ & ${\rm Gain}$  \\
\hline
$1$   &  $0.1$         &   $1.4$ & $4.13$ \\ 
\hline
$10$        &  $-0.64$     &   $-2.04$ & $2.11$ \\
\hline
$100$        &   $11$      & $0.02$  & $13.03$ \\
\hline
\end{tabular}
\caption{The performance of the trained three-point BDF scheme on the linear ODE \eqref{eq:odel} for three different values of the constant $c$. The gain in the fourth column is the ratio of the (mean) error \eqref{eq:odel4} with the standard BDF2 method and the (mean) error with the trained scheme \eqref{eq:bdfg} with parameters $g_{2,3}^{\ast}$, on the test set.}
\protect \label{tab:1}
\end{table}

The gradient descent algorithm converges quite quickly (atmost $9$ steps) to a local minimum of the non-convex loss function. The minimizers $g^{\ast}_{2,3}$ are shown in Table \ref{tab:1} for three different values of $c=1,10$ and $100$, and indicate a significant difference between the optimized values and the initial value of $(0.5,0.5)$ (corresponding to the standard second-order BDF2 scheme).

\begin{figure}[htbp]
\centering
\includegraphics[width=12cm]{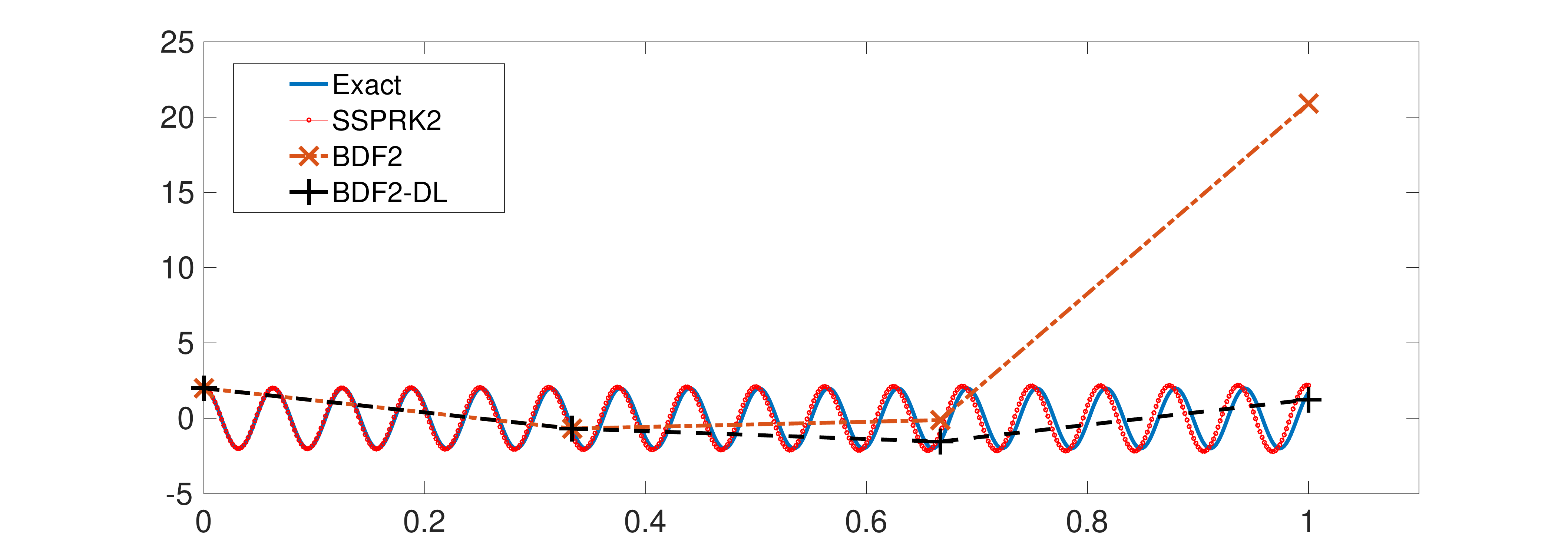}
\caption{Solutions $u$ of the linear ODE \eqref{eq:odel} (Y-axis) in time (X-axis), with $c=100$ in time period $[0,1]$ with an initial value from the \emph{test set}. We compare the exact solution with a solution computed with a second-order explicit Runge-Kutta method (SSP-RK2) with $1000$ time steps, the standard BDF2 method with $\Dt = \frac{1}{3}$ and the trained scheme (labelled as BDF2-DL (data learned)), also with $\Dt = \frac{1}{3}$. The trained scheme clearly outperforms the standard BDF2 method and more surprisingly, even the second-order accurate fine grid RK2 method. }
\label{fig:3}
\end{figure}

We test the trained scheme i.e, \eqref{eq:bdfg} with parameters $g^{\ast}_{2,3}$ on a \emph{test set}, chosen by randomly selecting $50$ points in the interval $[-5,5]$ as the initial data $u_0$ in \eqref{eq:odel}. The mean gains in error i.e the ratio of the error with the standard BDF2 scheme and the error with the trained (data learned) scheme, with respect to the underlying exact solution, are presented in table \ref{tab:1}. The gain in efficiency with the trained scheme is considerable for all the three values of $c$, rising to at least an order of magnitude for $c=100$. In this case, the value of $u$ computed with the trained scheme is remarkably close to the exact solution at times $T=2/3$ and $T=1$. Moreover, the trained scheme even seems to outperform an explicit second-order Runge-Kutta method (see \cite{LEV2}) with a fine grid of $\Dt = 0.001$, as observed in a plot of the approximate solutions on a particular realization of the test set in figure \ref{fig:3}. 
\subsection{A non-linear ODE.}
We consider a simple but non-linear population (saturation) model for the time evolution of the population density $u$, described by the ODE,
\begin{equation}
\label{eq:oden}
u^{\prime} = c u (1-u), \quad u(0) = u_0 \geq 0.
\end{equation}
It is straightforward to check that the exact solution of \eqref{eq:oden} is given by,
\begin{equation}
\label{eq:oden1}
u(t) = \frac{u_0}{u_0 + (1-u_0)e^{-ct}}.
\end{equation}
Hence, any non-negative initial condition converges to a saturation value (stable equilibrium) at $u=1$, at a time scale dictated by the constant $c$. We are interested in computing the solution in the time interval $[0,1]$. 

For step $1$ of Algorithm \ref{alg:1}, we again choose the generalized three-point BDF scheme \eqref{eq:bdfg} and a \emph{coarse} grid with time step $\Dt = 0.5$. To generate the training set, we choose initial data $\{u_0^i\}_{1\leq i \leq I}$ with $I=10$, uniformly over the interval $[0,2]$ and set $U_1^i$ as the exact solution \eqref{eq:oden1} at time $t= \Dt$. On this training set, we compute the approximate solution $U_2$ of \eqref{eq:oden} by the generalized three point BDF method (with a single trainable parameter $g_2$) at time $T=2\Dt=1$. The exact solution $U^{2,i}_{{\rm ref}}$ is calculated by \eqref{eq:oden1} at time $T=1$ with the initial data $u_0^i$ to define the \emph{loss function},
 \begin{equation}
\label{eq:oden2}
E_1(g_2) : = \sum\limits_{i=1}^I |U^{2,i}(g_2) - U^{2,i}_{{\rm ref}} |_1.
\end{equation}
Compared to the previous example of a linear ODE, we choose the $l^1$ norm as the loss function in this nonlinear example. 

\begin{table}[htbp]
\centering
\begin{tabular}{|c|c|c|}
\hline
$c$ & $g^{\ast}_2$ &  ${\rm Gain}$  \\
\hline
$0.2$   &  $0.41$         &   $1.5$ \\ 
\hline
$1$        &  $0.16$     &   $3.02$  \\
\hline
$5$        &   $0.03$      & $10.54$  \\
\hline
\end{tabular}
\caption{The performance of the trained three-point BDF scheme \eqref{eq:bdfg} on the non-linear ODE \eqref{eq:oden} on three different values of the constant $c$. The gain in the third column is the ratio of the (mean) error \eqref{eq:oden2} with the standard BDF2 method and the (mean) error with the trained scheme \eqref{eq:bdfg} with parameter $g_{2}^{\ast}$ on the test set.}
\protect \label{tab:2}
\end{table}

\begin{figure}[htbp]
\centering
\includegraphics[width=0.45\linewidth]{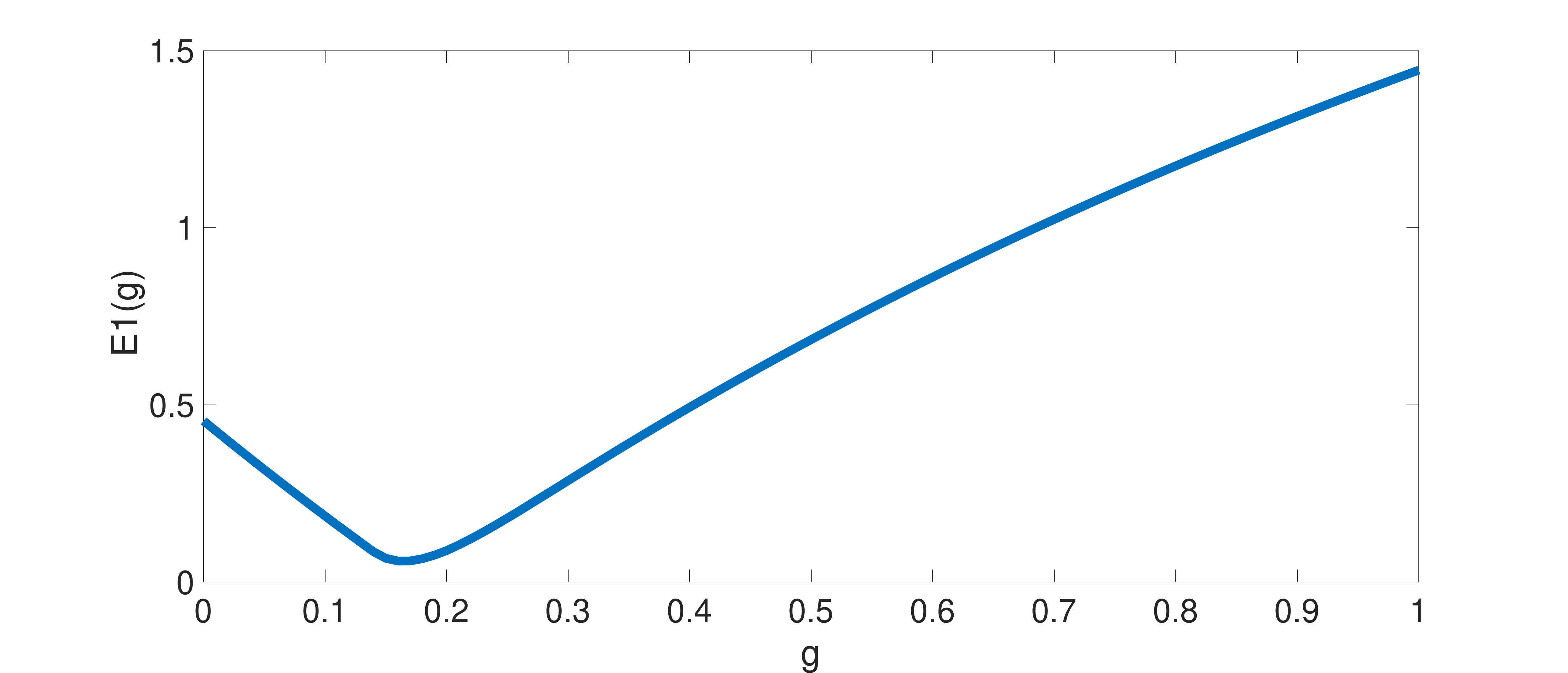} \includegraphics[width=0.51\linewidth]{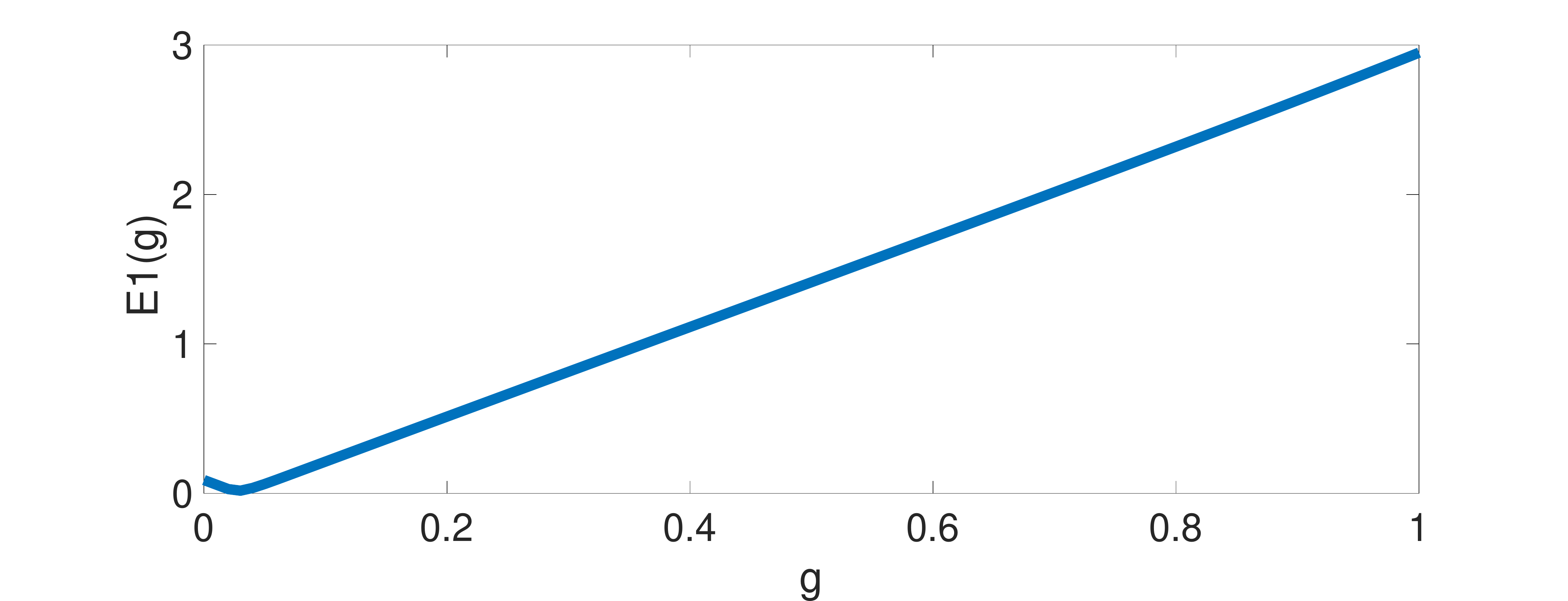}
\caption{Error $E_1(g)$ \eqref{eq:oden2} (Y-axis) vs. the scheme parameter $g$ (X-axis) in the generalized BDF scheme \eqref{eq:bdfg} for $\Dt = 0.5$ and two different values of $c$. Let $c=1$, the minimum in the error is achieved at
$g = 0.16$ Right: $c=5$: The minimum is at $g=0.03$.}
\label{fig:4}
\end{figure}

The loss function for two different values of $c=1$ and $c=5$ is shown in figure \ref{fig:4}. In all cases that we tested, the loss function is convex and is readily minimized by a straightforward steepest descent algorithm. The resulting optimal parameter $g_2^{\ast}$ for three different values of $c$ is shown in table \ref{tab:2}. As seen in table \ref{tab:2} and figure \ref{fig:4}, the optimal value $g^{\ast}$ is very different from $g=0.5$. 

The \emph{test set} is constructed by randomly choosing $50$ points from the interval $[0,5]$ as the initial data and approximating \eqref{eq:oden} with the trained scheme, i.e \eqref{eq:bdfg} with parameter $g_2^{\ast}$. The corresponding error with respect to the exact solution is calculated and the \emph{gain}, defined as before, is shown in table \ref{tab:2} . We observe a consistent gain in computational efficiency with the trained scheme that is approximately one order of magnitude for $c=5$. Thus, the machine learning algorithm \ref{alg:1} performs well in this nonlinear example and the gains in efficiency are similar to the linear problem even though a different loss function was used. 

\section{Heat equation}
\label{sec:5}
As the first example for PDEs, we consider the heat equation in one space dimension,
\begin{equation}
\label{eq:ht}
\begin{aligned}
u_t &= c u_{xx}, \quad (x,t) \in (0,1) \times (0,T), \\
u(x,0) &= u_0(x), \quad x \in (0,1), \\
u(0,t) &= u(1,t) = 0, \quad t \in (0,T).
\end{aligned}
\end{equation}
Here, $u$ is the temperature and $0 < c \in \R$ is a diffusion coefficient. 
\subsection{Numerical scheme.}
\label{sec:htnum}
We discretize the interval $[0,1]$ uniformly with a grid size $\Dx$ and label the resulting points as $x_j = j \Dx$ for $0\leq j \leq J+1$, with $\Dx = \frac{1}{J+1}$. The time interval $[0,T]$ is discretized uniformly with a time step $\Dt$ and the time points are labeled as $t^n = n \Dt$, with $0 \leq n \leq N$ and $\Dt = T/N$. We approximate the heat equation \eqref{eq:ht} by evolving $U^n_j \approx u(x_j,t^n)$, with the following generalized (or weighted) five-point finite difference scheme,
 \begin{equation}
 \label{eq:hts1}
 \begin{aligned}
 \frac{U^{n+1}_j - U^n_j}{\Dt} &= \frac{c(1-g^n)}{\Dx^2}\left(b^n_{-2} U^{n+1}_{j-2} + b^n_{-1}U^{n+1}_{j-1} + b^n_0 U^{n+1}_j + b^n_{1} U^{n+1}_{j+1} + b^n_{2} U^{n+1}_{j+2} \right) \\
                                               &+ \frac{cg^n}{\Dx^2} \left(b^n_{-2} U^{n}_{j-2} + b^n_{-1}U^{n}_{j-1} + b^n_0 U^{n}_j + b^n_{1} U^{n}_{j+1} + b^n_{2} U^{n}_{j+2} \right), \quad \forall n, 2 \leq j \leq J-1.
                                               \end{aligned}
                                               \end{equation} 
                                               
The update formulas for the points $U^n_{1}$ and $U^n_{J}$  is computed by setting the Dirichlet boundary conditions $U^{n}_{-1,0,J+1,J+2} \equiv 0$, for all $n$.

By using Taylor expansions, one can readily prove the following lemma,
\begin{lemma}
\label{lem:1}
For all $n$ and any $g^n \in \R$, the finite difference scheme \eqref{eq:hts1} is a consistent and first-order accurate discretization of the one-dimensional heat equation \eqref{eq:ht} if and only if the coefficients $b^n_k$, for $-2 \leq k \leq 2$, satisfy the following algebraic conditions,
\begin{equation}
\label{eq:hts2}
\begin{aligned}
b^n_2 + b^n_1 + b^n_0 + b^n_{-1} + b^n_{-2} &= 0, \\
2b^n_{2} + b^n_1 - b^n_{-1} - 2b^n_{-2} &= 0, \\
2b^n_2 +\frac{b^n_1}{2} + \frac{b^n_{-1}}{2} + 2b^n_{-2} &=0.
\end{aligned}
\end{equation}
\end{lemma}                                               
Here, consistency and accuracy are defined in terms of the local truncation error \cite{LEV2}. As \eqref{eq:hts2} has three equations containing five unknowns, we can eliminate three of them in terms of  $b^n_{-2,-1}$ to obtain,
\begin{equation}
\label{eq:hts3}
b^n_0 = 1-3b^n_{-1}-6b^n{-2}, \quad b^n_1 = 3b_{-1}^n + 8b_{-2}^n-2, \quad b_2^n =  1-b_{-1}^n-3b_{-2}^n.
\end{equation}

Hence, per time level, the scheme \eqref{eq:hts1} contains three undetermined parameters $g^n,b^n_{-1}$ and $b^n_{-2}$. These parameters will be determined by the training process, i.e, Step $2$ of Algorithm \ref{alg:1}.
\begin{remark}
Lemma \ref{lem:1} provides sufficient conditions for consistency of the finite difference scheme \eqref{eq:hts1}. Moreover, this scheme is \emph{conservative} i.e, $\sum_j U^{n+1}_j = \sum_j U^n_j$.  We can also obtain stability, for instance energy ($L^2$) stability or discrete maximum principles. These require additional constraints on the parameters and may constrain the training process further. We do not consider this aspect in the following. 
\end{remark}
Although the form \eqref{eq:hts1} of a two time-level, five point finite difference scheme is non-standard, it embeds several well-known finite difference approximations, namely
\begin{itemize}
\item [Scheme S1:] Backward Euler in time and second-order accurate in space by setting $g^n = 0, b^n_{-2} =0, b^n_{-1} = 1$ for all $n$.
\item [Scheme S2:] Crank-Nicolson in time and second-order accurate in space by setting $g^n = 0.5, b^n_{-2} =0, b^n_{-1} = 1$ for all $n$. 
\item [Scheme S3:] Backward Euler in time and fourth-order accurate in space by setting $g^n = 0, b^n_{-2} =-\frac{1}{12}, b^n_{-1} = \frac{4}{3}$ for all $n$. 
\item [Scheme S4:] Crank-Nicolson in time and fourth-order accurate in space by setting $g^n = 0.5, b^n_{-2} =-\frac{1}{12}, b^n_{-1} = \frac{4}{3}$ for all $n$.
\end{itemize}                                               
One can also set $g^n=1$ to recover an explicit forward Euler time discretization. However, we focus on implicit time stepping methods in order to avoid the constraint of  the severe CFL restriction for explicit time discretizations of the heat equation. 
\subsection{Training and Results on the test set.}
We approximate the solution $u$ of the heat equation with scheme \eqref{eq:hts1} at time $T=0.05$, for three different values of the diffusion coefficient $c$, namely $c=0.1,1,10$ ranging from slow to fast diffusion. All the experiments will be performed on a spatial grid with mesh size $\Dx = \frac{1}{10}$ i.e, with $10$ mesh points. Moreover, the time grid will based on a \emph{single} very large time step of $\Dt = 0.05$. 

We focus on varying the initial datum $u_0$ in \eqref{eq:ht} to generate the \emph{training set}. However, in contrast to ODEs, the initial datum $u_0$ for a PDE lies in an infinite dimensional function space, for instance $u_0 \in L^2((0,1))$. Given the challenge of approximating the resulting data to solution operator, in infinite dimensions, we focus on particular classes of initial datum, defined in terms of (finite dimensional) parameters. Motivated by applications in uncertainty quantification \cite{UQbook} and reduced order modeling \cite{RMbook}, we concentrate on the following specific parametric random initial data,

\subsubsection{Smooth data.}
\label{sec:kl1}
We consider the following $L$-term \emph{Karhunen-Loeve expansion},
\begin{equation}
\label{eq:kl1}
u_0(x,\omega) = \sum\limits_{l=1}^L \lambda_l Y_l(\omega) \sin(l\pi x),
\end{equation}
with $L=3$, $\lambda_l = \frac{1}{2^{l-1}}$ and the random numbers $Y_l(\omega)$ chosen from a uniform distribution on $[0,1]$. The training set is chosen by selecting (at random) $I$ draws of the random variables $Y_l^i(\omega)$ with $1\leq i \leq I=20$. We compute a reference solution, for the resulting initial data $u_0^i$, with an explicit forward Euler time stepping and standard second-order spatial finite difference discretization \cite{LEV2} on a very fine grid of $1000$ mesh points and a time step, chosen to satisfy the standard CFL requirement for the heat equation. This fine grid solution is projected onto the underlying coarse grid by sampling this solution at points $x_j$ and at final time $T=\Dt = 0.05$, $1\leq j \leq J$. We denote this reference solution as $U^{n,i}_{j,{\rm ref}}$. 

The \emph{loss function} is defined as the $L^2$ error,
\begin{equation}
\label{eq:htlf}
E_2\left(g^1,b^1_{-2},b^1_{-1}\right) := \frac{\Dx}{2} \sum\limits_{i=1}^{I} \sum\limits_{j=1}^J |U^{1,i}_j - U^{1,i}_{j,{\rm ref}}|^2.
\end{equation} 

The loss function \eqref{eq:htlf} is minimized using a simplified version of the stochastic gradient algorithm \cite{SG} with a batch size of $4$, initialized with the starting values of $g^1 = 0.5$, $b^1_{-2} = 0$ and $b^1_{-1} = 1$, corresponding to the overall second-order Crank-Nicolson type scheme S2.  We denote the (approximate) minimizers as $\{g^{1,\ast},b^{1,\ast}_{-2},b^{1,\ast}_{-1}\}$ and the trained (\emph{data learned}) scheme is the finite difference \eqref{eq:hts1} with these parameters. 

The training (for three different cases of the diffusion coefficient $c$ considered here) resulted in (local) minimizers shown in table \ref{tab:3}. We observe that in most cases, the trained scheme is very different from any standard scheme. A \emph{test set} is generated by choosing $100$ random values of $Y_l$, $l=1,2,3$ in \eqref{eq:kl1}. Care is taken to exclude repetition of values from the training set and the corresponding reference solution is computed, analogously to the training set. 

\begin{figure}[htbp]
\centering
\includegraphics[width=12cm]{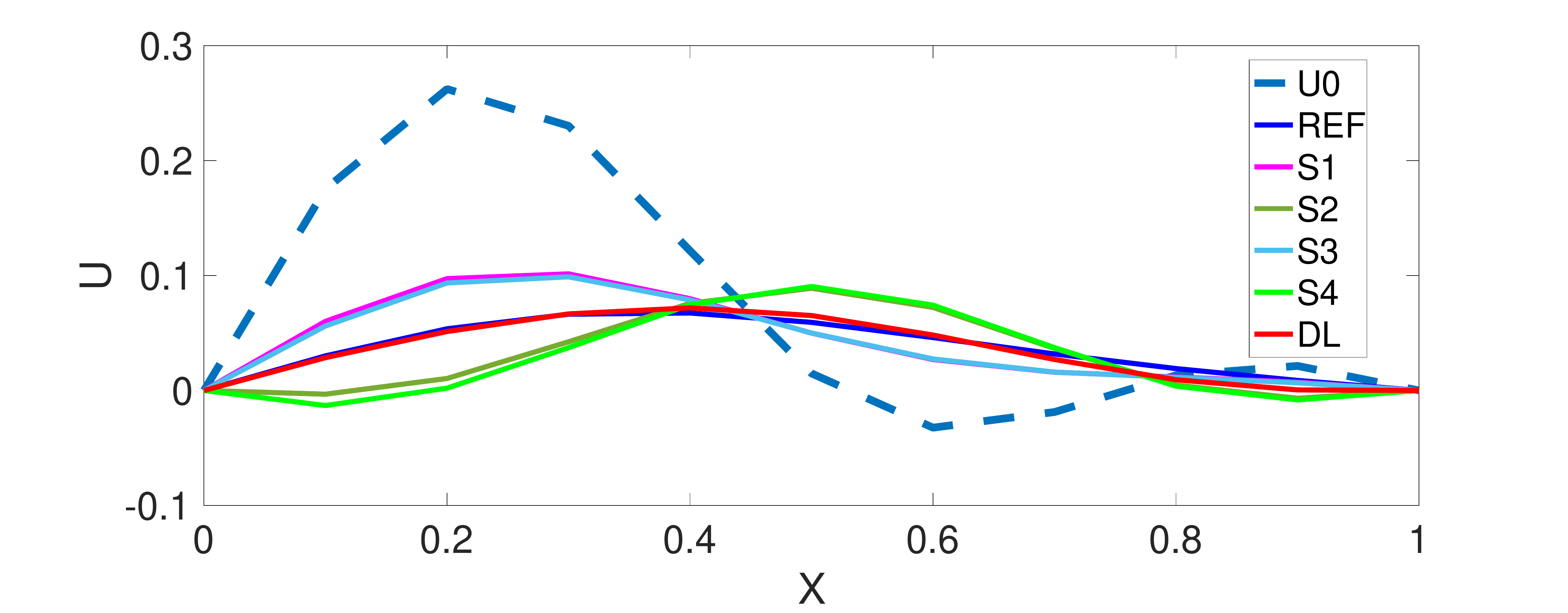}
\caption{Solutions of the heat equation \eqref{eq:ht} with $c=1$ at time $T=0.05$ with an initial value from the \emph{test set} for smooth data i.e \eqref{eq:kl1} with a particular realization of $Y_{1,2,3}(\omega)$. We show the initial data (U0), reference solution (REF), solutions approximated with schemes S1 (second-order in space+Backward Euler), S2 (second-order in space+Crank-Nicolson), S3 (fourth-order in space + backward Euler), S4 (fourth-order in space + Crank Nicolson) with the solutions by computed \emph{Data learned} (DL) trained scheme, \eqref{eq:hts1} with parameters given in table \ref{tab:3}, on a grid with one time step and $10$ mesh points. Observe the considerable gain in accuracy with the DL scheme over all the other schemes. }
\label{fig:5}
\end{figure}

Summarizing these results, we first observe that there is no clear winner among the standard schemes on the test set . For slow to moderate values of the diffusion coefficient, the scheme S2 i.e, Crank-Nicolson in time and second-order in space scores over the other three schemes whereas for a large value of the diffusion coefficient, the scheme S1 and S3 clearly perform the best. On the other hand, there is a large gain with the trained scheme compared to all the standard schemes. The gain of approximately $4$ is most modest for the $c=1$ value of the diffusion coefficient. On the other hand, there is clearly a gain of a factor of at least $10$ or $20$ for the extreme values of the diffusion coefficient. This gain in accuracy comes at no additional online cost and justifies the efficacy of the proposed machine learning algorithm.  This significant gain in performance with the trained scheme, for one particular instance of the test set, is also displayed in figure \ref{fig:5}.
\begin{table}[htbp]
\centering
\begin{tabular}{|c|c|c|c|c|c|c|c|}
\hline
$c$ & $g^{1,\ast}$ & $b^{1,\ast}_{-2}$ & $b^{1,\ast}_{-1}$  & Gain-1 & Gain-2 & Gain-3 & Gain-4 \\
\hline
$0.1$   &  $0.64$         &   $0$ & $1$ & $77.31$ & $19.58$ & $49.6$ & $38.97$ \\ 
\hline
$1$        &  $0.3$     &   $0$ & $1.12$ & $3.97$ & $3.21$ & $3.59$ & $3.98$ \\
\hline
$10$        &   $0.04$      & $1.2$  & $2.43$ & $11.91$ & $46.97$ & $11.49$ & $48.17$ \\
\hline
\end{tabular}
\caption{The performance of the trained finite difference scheme \eqref{eq:hts1} on the heat equation \eqref{eq:ht} for three different values of the constant $c$ and for the smooth initial data \eqref{eq:kl1}. The gain-m is the ratio of the (mean) error with the standard S-m scheme and the (mean) error with the trained scheme i.e, \eqref{eq:hts1} with parameters $g^{1,\ast},b^{1,\ast}_{-2,-1}$, on the test set.}
\protect \label{tab:3}
\end{table}

\subsubsection{Rough data}
Next, we consider the following discontinuous random initial data,
\begin{equation}
\label{eq:kl2}
u_0(x,\omega) = \begin{cases} 
                           1 + \epsilon Y_1(\omega), &{\rm if} \quad \frac{1}{3} + \epsilon Y_2(\omega) < x < \frac{2}{3} + \epsilon Y_3(\omega), \\
                           0, &{\rm otherwise},
                           \end{cases}
\end{equation}                           
Here, $\epsilon = 0.2$ and $Y_{1,2,3}$ are chosen randomly from a uniform distribution on $[-1,1]$. In other words, the initial data \eqref{eq:kl2} represents a step function with two discontinuities where the amplitude of the jump at the discontinuity and the location of both jumps are random. The underlying coarse grid is the same as for the smooth case. The training and test sets are generated in a manner, identical to the smooth case and the loss function \eqref{eq:htlf} is minimized similarly.

The training (for three different cases of the diffusion coefficient $c$ considered here) resulted in (local) minimizers shown in table \ref{tab:4}. We report that the training process converged very slowly for the $c=10$ value in the case of this rough initial data. This is reflected in the values of $20$ for both spatial weights as we terminated the iterations in the gradient descent method at this stage. One can provide a heuristic explanation for the values of the minimizers in this case. Recall that $c=10$ implies a very large amount of diffusion in the solution, such that the solution is almost zero at $T=0.05$. The value of $g^1 = 0$ corresponds to the most diffusive backward Euler method and similarly very high values for $b^1_{-2,-1}$ also imply a large amount of diffusion and drive the approximate solution (computed by \eqref{eq:hts1}) to zero. 

The gains with the trained scheme, over the four standard schemes, are shown in table \ref{tab:4} and indicate a very large gain over the best performing of the standard schemes, amounting to a factor of approximately $50$ for the case of $c=10$.

\begin{table}[htbp]
\centering
\begin{tabular}{|c|c|c|c|c|c|c|c|}
\hline
$c$ & $g^{1,\ast}$ & $b^{1,\ast}_{-2}$ & $b^{1,\ast}_{-1}$  & Gain-1 & Gain-2 & Gain-3 & Gain-4 \\
\hline
$0.1$   &  $0.24$         &   $-0.26$ & $2.15$ & $12.7$ & $3.16$ & $10.56$ & $3.78$ \\ 
\hline
$1$        &  $0.14$     &   $-0.38$ & $2.53$ & $2.09$ & $6.03$ & $1.98$ & $6.38$ \\
\hline
$10$        &   $0$      & $20$  & $20$ & $46.02$ & $228.93$ & $44.72$ & $231.32$ \\
\hline
\end{tabular}
\caption{The performance of the trained finite difference scheme \eqref{eq:hts1} on the heat equation \eqref{eq:ht} for three different values of the constant $c$ and for the rough initial data \eqref{eq:kl2}. The gain-m is the ratio of the (mean) error with the standard S-m scheme and the (mean) error with the trained scheme i.e, \eqref{eq:hts1} with parameters $g^{1,\ast},b^{1,\ast}_{-2,-1}$, on the test set.}
\protect \label{tab:4}
\end{table}
\section{Linear advection equation}
\label{sec:4}
The linear advection equation is considered as a prototype for the design and analysis of efficient numerical methods for hyperbolic equations. In one space dimension, it is given by
\begin{equation}
\label{eq:lt}
\begin{aligned}
u_t + c u_x &=0, \quad (x,t) \in (0,1) \times (0,T), \\
u(x,0) &= u_0(x), \quad x \in (0,1).
\end{aligned}
\end{equation}
For definiteness, we assume that $0 \leq c \in \R$ and supplement \eqref{eq:lt} with periodic boundary conditions. 
\subsection{Numerical scheme.}
We discretize the computational domain $[0,1] \times [0,T]$ as in section \ref{sec:htnum} and  use the following three-point finite difference scheme to approximate the linear advection equation \eqref{eq:lt} by evolving $U^n_j \approx u(x_j,t^n)$ with 
 \begin{equation}
 \label{eq:lts1}
 \begin{aligned}
 \frac{U^{n+1}_j - U^n_j}{\Dt} &= \frac{c(1-g^n)}{\Dx}\left(b^n_{-1}U^{n+1}_{j-1} + b^n_0 U^{n+1}_j + b^n_{1} U^{n+1}_{j+1} \right) \\
                                               &+ \frac{cg^n}{\Dx} \left(b^n_{-1}U^{n}_{j-1} + b^n_0 U^{n}_j + b^n_{1} U^{n}_{j+1}\right), \quad \forall n, 1 \leq j \leq J.
                                               \end{aligned}
                                               \end{equation} 
The update formulas for the points $U^n_{1}$ and $U^n_{J}$  are computed by using the periodic boundary conditions.
By using Taylor expansions, one can readily prove the following lemma,
\begin{lemma}
\label{lem:2}
For all $n$ and any $g^n \in \R$, the finite difference scheme \eqref{eq:lts1} is a consistent and first-order accurate discretization of the linear advection equation \eqref{eq:lt} if and only if the coefficients $b^n_k$, for $-1 \leq k \leq 1$, satisfy the following algebraic conditions,
\begin{equation}
\label{eq:lts2}
b^n_1 + b^n_0 + b^n_{-1} = 0, \quad b^n_1 - b^n_{-1} = 1.
\end{equation}
\end{lemma}                                               
We can eliminate two parameters in the system \eqref{eq:lts2} in terms of the undetermined parameter $b^n_{-1}$ to obtain,
\begin{equation}
\label{eq:lts3}
b^n_0 = -1-2b_{-1}, \quad b^n_1 = 1+b^n_{-1}
\end{equation}

Hence, per time level, the scheme \eqref{eq:lts1} contains two undetermined parameters $g^n$ and $b^n_{-1}$. The scheme \eqref{eq:lts1} with constraints \eqref{eq:lts2} is consistent as well as \emph{conservative} i.e, $\sum_j U^{n+1}_j = \sum_j U^n_j$. Additional constraints on the parameters are needed to impose stability conditions such as discrete $L^2$ (energy) stability or a discrete maximum principle.

The generalized form \eqref{eq:lts1} embeds several well-known finite difference approximations namely,
\begin{itemize}
\item [Scheme S1:] Backward Euler in time and upwind in space by setting $g^n = 0, b^n_{-1} =0$
\item [Scheme S2:] Crank-Nicolson in time and upwind in space by setting $g^n = 0.5, b^n_{-1} =0$  
\item [Scheme S3:] Backward Euler in time and central in space by setting $g^n = 0, b^n_{-1} = 0.5$. 
\item [Scheme S4:] Crank-Nicolson in time and central accurate in space by setting $g^n = 0.5, b^n_{-1} = 0.5$.
\end{itemize}        
As for the heat equation, we consider only implicit (in time) schemes as they do not require a restriction on the time step $\Dt$. We approximate the solution $u$ of the linear advection equation with scheme \eqref{eq:lts1} at time $T=0.5$ and consider two different values of the wave speed $c$, namely $c=0.5$ and $c=2$  All the experiments will be performed on a very coarse spatial grid with mesh size $\Dx = \frac{1}{10}$ i.e, $10$ mesh points, and a single, very large, time step of $\Dt = 0.5$. 

\begin{table}[htbp]
\centering
\begin{tabular}{|c|c|c|c|}
\hline
$c$ & $g^{1,\ast}$ & $b^{1,\ast}_{-1}$  & Gain \\
\hline
$0.5$   &  $0.2$         &   $-2$ & $3.72$ \\ 
\hline
$2.0$   &  $-20$         &   $0$ & $96.51$ \\ 
\hline
\end{tabular}
\caption{The performance of the trained finite difference scheme \eqref{eq:lts1} on the linear advection equation \eqref{eq:lt} for two different values of the wave speed $c$ and for initial data \eqref{eq:kl1}. The gain is the ratio of the (mean) error with the best performing of the (S1,S2,S3,S4) schemes and the (mean) error with the trained scheme i.e, \eqref{eq:lts1} with parameters $g^{1,\ast},b^{1,\ast}_{-1}$, on the test set.}
\protect \label{tab:5}
\end{table}

\subsection{Training and Results on the test set.}
To generate the training and test sets, we use an identical set up as described for the heat equation in section \ref{sec:kl1}. In particular, both the training and test sets are generated from the three term Karhunen-Loeve expansion \eqref{eq:kl1}. We compute a reference solution on a fine mesh of $1000$ points, with an explicit forward Euler time stepping and standard upwind finite difference discretization \cite{LEV2}. The time step is chosen to satisfy the standard CFL requirement for the advection equation. This fine grid solution is projected onto the underlying coarse grid by sampling this solution at points $x_j$ and at final time $T=\Dt = 0.5$, $1\leq j \leq J$, to generate the reference solutions. The loss function \eqref{eq:htlf}  is minimized with a simplified version of the stochastic gradient algorithm \cite{SG} with a batch size of $4$. We initialize the gradient descent method with the starting values of $g^1 = 0.5$, $b^1_{-1} = 0$, corresponding to the Crank-Nicolson in time, upwind in space, scheme S2,  and denote the (approximate) minimizers as $\{g^{1,\ast},b^{1,\ast}_{-1}\}$. The minimizers, for different values of $c$, are shown in table \ref{tab:5}. We observe that for $c=0.5$, the computed minimizers deviate greatly from any standard scheme. However, this contrast is much more pronounced in the $c=2$ case as there was very slow convergence of the stochastic gradient method and it was terminated at $g^{1,\ast} = -20$, indicating that there is a path along which the loss function (very slowly) approaches a value of zero.

\begin{figure}[htbp]
\centering
\includegraphics[width=0.45\linewidth]{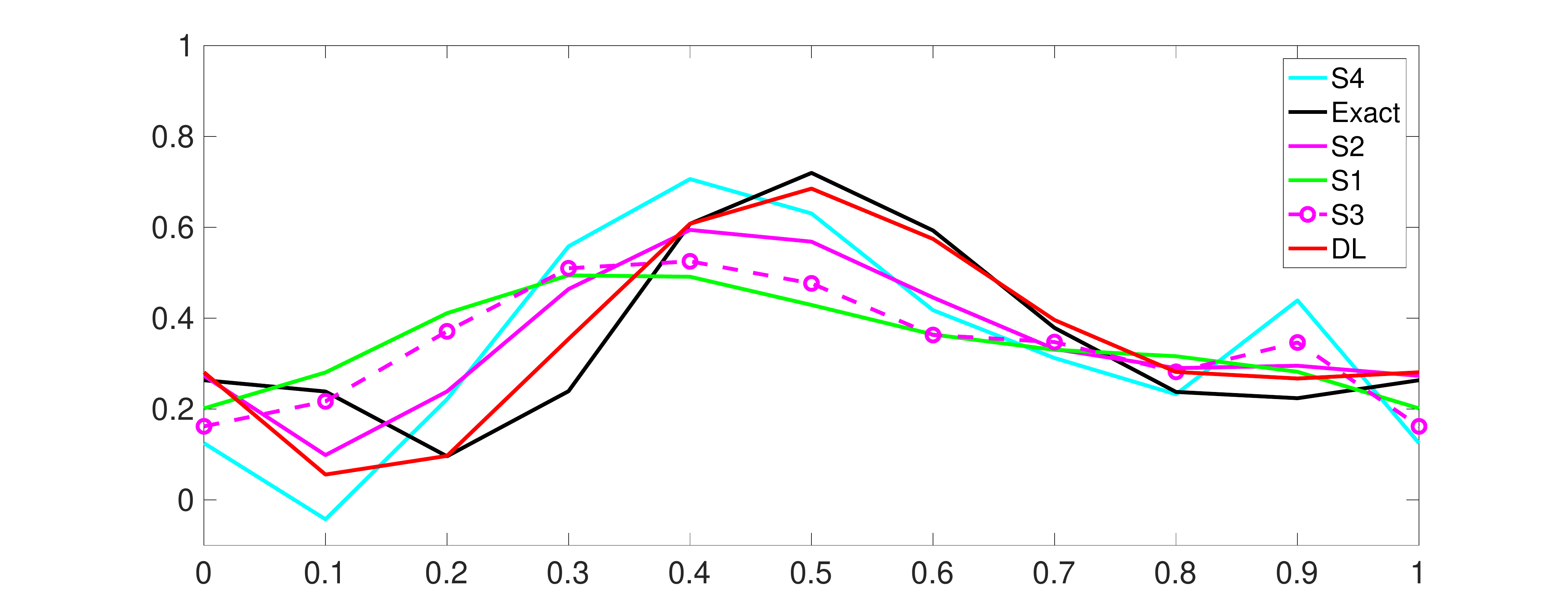} \includegraphics[width=0.5\linewidth]{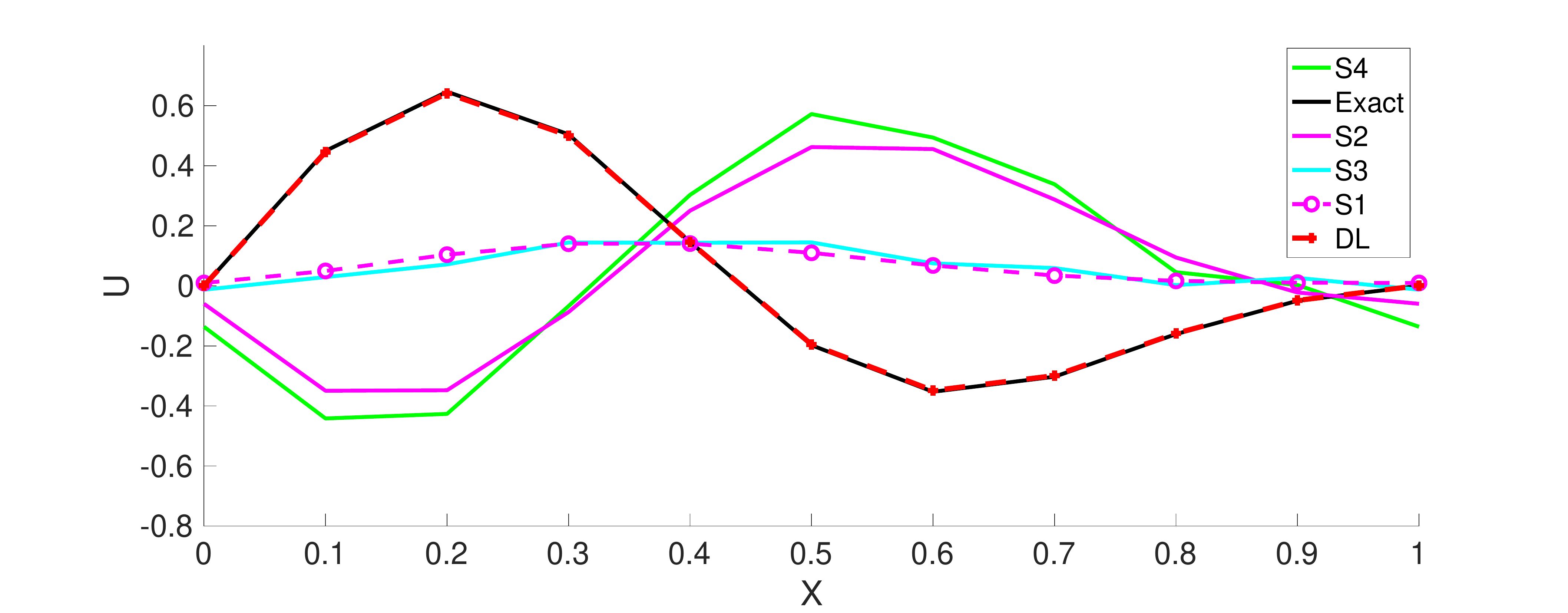}
\caption{Comparison of five different schemes for the linear advection equation \eqref{eq:lt}, namely Backward Euler in time and upwind in space (S1), Crank-Nicolson in time and upwind in space (S2), Backward Euler in time and central in space (S3), Crank-Nicolson in time and central in space (S4) and the trained scheme (DL) i.e, \eqref{eq:lts1} with weights given in table \ref{tab:5} on grid of $10$ mesh points and with one time step. Comparison with the exact solution computed with the upwind scheme with forward Euler time stepping on a fine grid of $500$ points, and on a randomly chosen initial data from \eqref{eq:kl1}, for two different wave speeds. Left $c=0.5$, Right $c=2$.}
\label{fig:6}
\end{figure}

In order to compare with standard schemes, we define a ${\rm Gain}$ as the ratio of the (mean) error on the test set with the best performing of the four schemes $S1,S2,S3,S4$ (the one with the least mean error) and the trained scheme. For $c=0.5$, the scheme S2 is the best performing scheme and for $c=2$, the scheme $S1$ is the best performing scheme. The computed gain is shown in table \ref{tab:5}.  For further comparison, we plot a single randomly chosen realization of the test data for both $c=0.5$ and $c=2$ in figure \ref{fig:6}. 

As shown in table \ref{tab:5}, for the case of $c=0.5$, the trained scheme provides a gain of $3.72$ over the best performing of the standard schemes (the scheme S2). The gains with respect to the backward Euler time stepping schemes are larger. This is also shown in figure \ref{fig:6} (left), where we observe that the trained scheme is significantly more accurate than standard schemes. 

However, the gains with the trained scheme are enormous in the case of $c=2$, amounting to a gain of almost two orders of magnitude vis a vis the best performing of the standard schemes, see table \ref{tab:5} and figure \ref{fig:6} (right). A heuristic explanation for this observation goes as follows: recall that the exact solution coincides with the initial data in this case. Thus in the limit of $g^n \rightarrow -\infty$, we can see from \eqref{eq:lts1} that $U^{n+1} \approx U^n$ and we are very close to the initial data. It appears that the machine learning algorithm \emph{learns} this fact, when shown training data, and provides this remarkable gain in accuracy in this special case. 
\section{Burgers' equation.}
\label{sec:6}
The Burgers' equation given by
\begin{equation}
\label{eq:bg}
\begin{aligned}
u_t + \left(\frac{u^2}{2}\right)_x &= 0, \quad (x,t) \in (0,1) \times (0,T), \\
\end{aligned}
\end{equation}
is a prototypical example for nonlinear hyperbolic conservation laws,
\begin{equation}
\label{eq:cl}
\begin{aligned}
u_t + (f(u))_x &= 0, \\
u(x,0) &= u_0(x),
\end{aligned}
\end{equation}
These equations arise in a wide variety of applications and examples include the Euler equations of gas dynamics, the shallow water equations of oceanography and the MHD equations of plasma physics \cite{DAF1}. It is well-known that solutions of \eqref{eq:cl} develop finite time singularities in the form of \emph{shock waves}, when even the initial data is smooth. Thus, solutions of \eqref{eq:cl} are sought in the sense of distributions and additional \emph{entropy} conditions are imposed in order to recover uniqueness \cite{DAF1}.
\subsection{Numerical Scheme}
\label{sec:clnum}
There is a large body of literature on numerical methods for hyperbolic conservation laws and popular numerical methods include the conservative finite difference schemes and discontinuous Galerkin finite element methods \cite{GR1}. However, for the sake of simplicity, we consider the simplest \emph{first order finite volume scheme} in this section. 

We discretize the interval $[0,1]$ uniformly with a grid size $\Dx$ and label the resulting points as $x_j = j \Dx$ for $0\leq j \leq J+1$, with $\Dx = \frac{1}{J+1}$. Thus, the interval is partitioned into cells (control volumes),
$$
C_{j} := (x_{\jmhf},x_{\jphf}), \quad x_{\jphf} = \frac{x_{j} + x_{j+1}}{2}, ~ \forall j.
$$
The time interval $[0,T]$ is discretized uniformly with a time step $\Dt$ and the time levels are denoted as $t^n = n \Dt$. We approximate the cell averages of the solution of \eqref{eq:cl},
$$
U^n_j \approx \int_{C_j} u(x,t^n) dx,
$$
by writing the update formula,
\begin{equation}
\label{eq:bfv1}
\begin{aligned}
U^{n+1}_j &= U^n_j - \frac{\Dt}{\Dx}\left(F^n_{\jphf} - F^n_{\jmhf} \right), \\
U_j^n &=  \int_{C_j} u_0(x) dx
\end{aligned}
\end{equation}
Here, $F^n_{\jphf} = F(U^n_j,U^n_{j+1})$ is a numerical flux, consistent with the flux function $f$ in \eqref{eq:cl}. A popular choice is the so-called Local Lax-Friedrichs or \emph{Rusanov} flux given by,
\begin{equation}
\label{eq:rus1}
\begin{aligned}
F(U^n_j,U^n_{j+1}) &= \frac{1}{2}\left(f(U^n_j) + f(U^n_{j+1})\right) - \frac{1}{2}\max(|f^{\prime}(U^n_j)|,|f^{\prime}(U^n_{j+1})|)(U^n_{j+1} - U^n_j).
                                \end{aligned}
\end{equation}
Thus, the numerical diffusion is weighted by a local \emph{wave speed}. It is well-known that the resulting scheme \eqref{eq:bfv1} with flux \eqref{eq:rus1} is \emph{conservative, consistent} and \emph{monotone} \cite{GR1}. Consequently, the solutions computed by the scheme converge to an entropy solution of \eqref{eq:cl}.

We cast the finite volume scheme \eqref{eq:bfv1} in our machine learning framework by generalizing the flux \eqref{eq:rus1} to
\begin{equation}
\label{eq:rusdl1}
\begin{aligned}
F(U^n_j,U^n_{j+1}) &= \frac{1}{2}\left(f(U^n_j) + f(U^n_{j+1})\right) - \left(w^n_{\jphf}\max(|f^{\prime}(U^n_j)|,|f^{\prime}(U^n_{j+1})|)(U^n_{j+1} - U^n_j)\right).
\end{aligned}
\end{equation}
Here $w^n_{\jphf} \in \R,~\forall j$, are weights corresponding to the scaling of the local wave speed. The resulting scheme is given by,
\begin{equation}
\label{eq:bfvl}
 \begin{aligned}
U^{n+1}_j &= U^n_j - \frac{\Dt}{2\Dx} \left(f(U^n_{j+1})-f(U^n_{j-1})\right) \\
&+\frac{\Dt}{\Dx} \left(w^n_{\jphf}\max(|f^{\prime}(U^n_j)|,|f^{\prime}(U^n_{j+1})|)(U^n_{j+1} - U^n_j) - w^n_{\jmhf}\max(|f^{\prime}(U^n_j)|,|f^{\prime}(U^n_{j-1})|)(U^n_{j} - U^n_{j-1})\right).
\end{aligned}
\end{equation}
The properties of this \emph{generalized finite volume scheme} \eqref{eq:bfvl} are summarized in the lemma below,
\begin{lemma}
\label{lem:3}
The solutions $U^n_j$ generated by the scheme \eqref{eq:bfvl} satisfy the following,
\begin{itemize}
\item [i.] For any $w^n_{\jphf} \in \R$, for all $j,n$, the scheme \eqref{eq:bfvl} is conservative and consistent. Moreover, it is (formally) first-order accurate. 
\item [ii.] Under the CFL condition,
\begin{equation}
\label{eq:cfl1}
\max_{j} |f^{\prime}(U^n_j)| \frac{\Dt}{\Dx} \leq 1,
\end{equation}
and under the condition $w^n_{\jphf} \geq 0.5$, for all $j,n$, the scheme \eqref{eq:bfvl} is monotone. Hence, the corresponding approximations converge to the entropy solution of the scalar conservation law \eqref{eq:cl}.
\end{itemize}
\end{lemma}
The proof of the above lemma is a straightforward adaptation of the results of \cite{GR1}. 

\begin{remark}
We remark that setting $w^n_{\jphf} \equiv 0.5$ for all $n,j$, yields the standard Rusanov scheme. The weights $w$ in \eqref{eq:bfvl} serve to modulate the numerical diffusion in the scheme. We can design an alternative
scheme by replacing the arithmetic averages of the fluxes in \eqref{eq:rusdl1} by an \emph{entropy conservative flux} and replacing the conservative variable $u$ in the numerical diffusion term of the flux \eqref{eq:rusdl1} by the entropy
variables \cite{TAD1,FMT4}. The resulting scheme can be proved to be \emph{entropy stable} for any $w^n_{\jphf} \geq 0$. 
\end{remark}

Motivated by machine learning architectures such as convolution neural networks \cite{DLbook}. we make a further simplification by \emph{pooling} values of the weights $w^n$ in longer \emph{windows} i.e by requiring that
\begin{equation}
\label{eq:wind}
w^n_{j+1/2} = w^n_{j_s + 1/2}, \forall |j-j_s| \leq m.
\end{equation}
Here $0\leq m \leq J$ is the window length and $j_s \subseteq \{1,\ldots, J\}$ is a subset of grid points on which we center the pooling windows. 
\subsection{Generalized Rusanov scheme as a neural network.}
\label{sec:rusnn}
\begin{figure}[htbp]
\centering
\includegraphics[width=10cm]{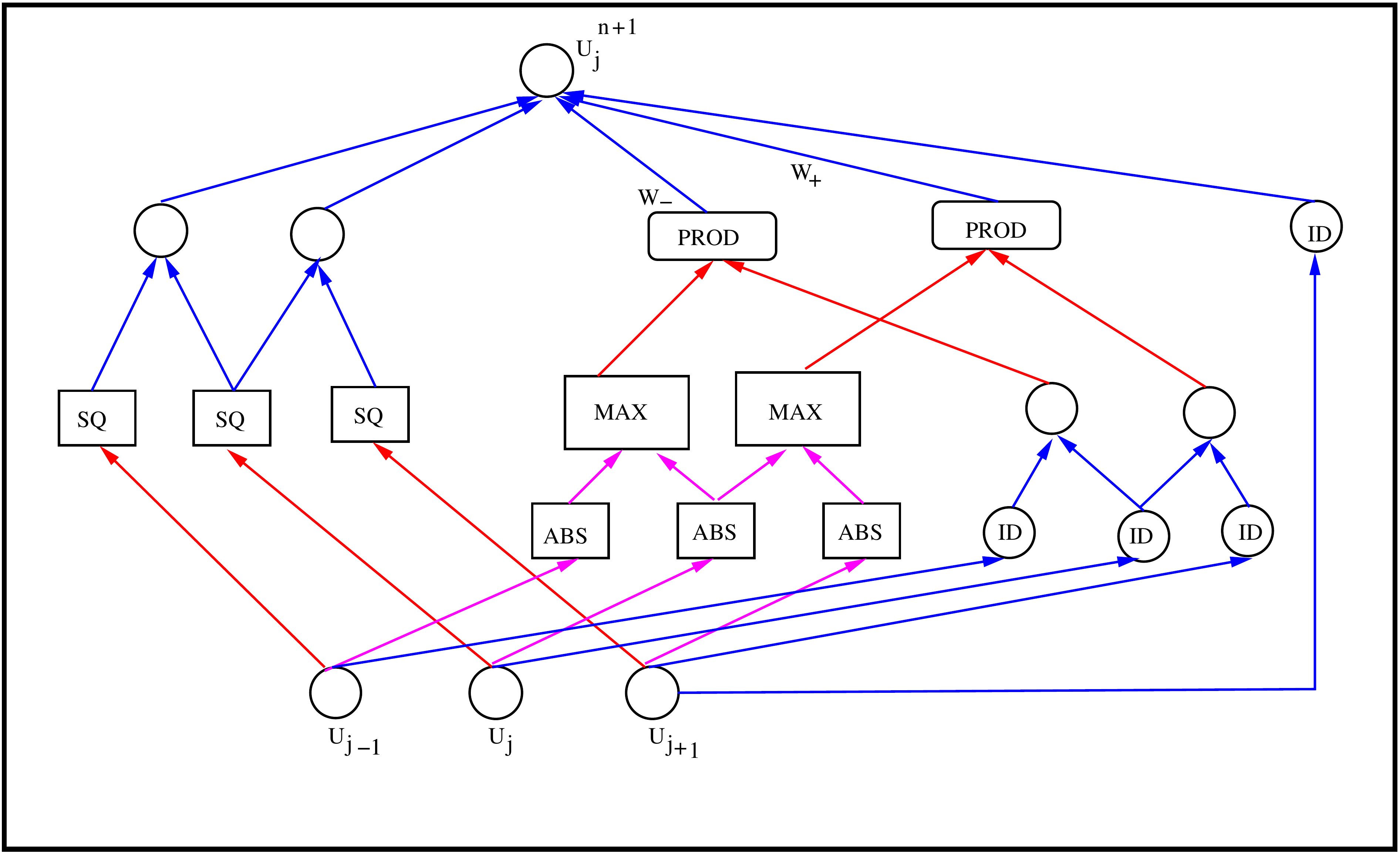}
\caption{Representation of a single time step of the Rusanov type scheme \eqref{eq:bfvl} for approximating the Burgers' equation \eqref{eq:bg} as a \emph{multi-layer neural network}. As inputs, we consider the three point stencil $\{U^n_{j-1}, U^n_j, U^n_{j+1}\}$ that are converted to the output $U^{n+1}_j$ through a series of operations. Linear operations are shown with blue arrows. We label SQ to represent the function $f(u) = 0.5u^2$, PROD to be product of two numbers, MAX as the maximum of two numbers and ABS to represent $|u|$. Both MAX and ABS are directly represented by a combination of linear and RELU steps and are shown by magenta arrows where as SQ and PROD are shown with red arrows. The only undetermined weights $W^{\pm}$ are those of $w_{j\pm 1/2}^n$ in \eqref{eq:bfvl}. The linear operation $Au = u$ is labeled as ID. Note that SQ and PROD can be very efficiently approximated by a rather narrow sub-network based on ReLU activation units \cite{YAR1}.} 
\label{fig:7}
\end{figure}

The generalized Rusanov scheme \eqref{eq:bfvl} can be represented as an artificial neural network as shown in figure \ref{fig:7}, where we focus on the specific example of the Burgers' equation \eqref{eq:bg}. Note that in figure \ref{fig:7}, the input neurons or units are components of the vector of unknowns $U^n = \{U^n_j\}$. The output at the end of each time step is the vector $U^{n+1}$. The network transforms the input into the output through layers of operations that consist of linear maps (connecting different neurons) and nonlinear operations. For this specific example, there are the following nonlinear operations,
\begin{itemize}
\item [ABS]: refers to $ABS(a) = |a| = \sigma(a) + \sigma(-a)$, with $\sigma$ being the ReLU activation function \eqref{eq:relu}.
\item [MAX]: refers to $MAX(a,b) = \max(a,b) = a+\sigma(b-a)$.
\item [SQ]; refers to $SQ(a) = 0.5a^2$.
\item [PROD]: refers to $PROD(a,b) = ab$.
\end{itemize}
Thus, ABS and MAX are directly expressed in the terms of a \emph{traditional} neural network, in the sense of \cite{DLbook}, with a ReLU activation. On the other hand, SQ and PROD are bespoke nonlinear operations for this particular example. Hence, the neural network in figure \ref{fig:7} is not a traditional neural network. However, it can be directly represented in terms of the so-called \emph{sum product networks} \cite{SP}. On the other hand, following recent papers such as \cite{YAR1}, one can approximate the square and product functions very efficiently in terms of neural networks with ReLU activations, see \cite{SZ1} figure 1. In particular, one can approximate the square and product maps up to accuracy $\delta$ by a neural network of width that is at most logarithmic in $\delta$. Therefore, the neural network underlying the scheme \eqref{eq:bfvl} is very readily approximated by a traditional ReLU based network architecture. Given several hidden layers per time step (see figure \ref{fig:7}) and multiple time steps, the scheme \eqref{eq:bfvl} is realized as a \emph{deep neural network}.

It is standard in deep learning to optimize the weights (corresponding to entries in all matrices) for the whole network during the training process. It is in this step that we differ from traditional machine learning and take a more conservative approach. Given that we wish to be consistent (and formally first-order accurate) for any of the weights that might crop up in the training process, we severely constrain the set of free parameters within the network to only the weights $w$ of the numerical viscosity coefficient in \eqref{eq:bfvl}. Even these weights are pooled in windows. Thus, at most we have two free (trainable) parameters in the sub-network shown in figure \ref{fig:7}. Hence, the training process operates on a   (considerably) constricted part of the deep neural network. 

\subsection{Training and results on the test set.}
For the numerical experiments, we will only consider the Burgers' equation \eqref{eq:bg} with periodic boundary conditions. We fix $\Dx = 0.1$ i.e, we discretize the interval $[0,1]$ into $10$ cells. Our final time is $T=0.1$ the time interval is discretized into \emph{two time steps} with time step size $\Dt = 0.05$. On this coarse spatial grid, such a large time step is consistent with the CFL condition \eqref{eq:cfl1}. Moreover, we \emph{pool} the weights of the numerical diffusion operator by a defining a pooling window of size $m=3$ in \eqref{eq:wind}. Hence, at each time step we need to specify $3$ weights namely $\{w^n_{1},w^n_{2},w^n_{3}\}$, correspond of the first, middle and last three of the interior interfaces. The weights on the boundary interfaces are determined from the periodic boundary conditions.
 
To generate the training set, we first consider the \emph{smooth} random initial data, specified by the Karhunen-Loeve expansion \eqref{eq:kl1}. As in the previous sections, the expansion is truncated at $L=3$ and we choose random variables $Y^i_{1,2,3}$, for $1 \leq i \leq I$ and $I=20$, uniformly from $[0,1]$. The corresponding initial data, labeled as $u_0^i$ is used to initialize the scheme \eqref{eq:bfvl} and to generate the updated solutions $U^{1,i}_j,U^{2,i}_j$ for all $j$. A reference solution is computed using the standard Rusanov scheme (i.e setting $w^n_{\jphf} \equiv 0.5$) on a fine mesh of $1000$ points and with the time step determined by the corresponding CFL number as in \eqref{eq:cfl1}. This fine grid solution at the two time levels $\Dt,2\Dt$ is projected to the underlying coarse grid by averaging over the coarse cells and is denoted as $U^{n,i}_{j,{\rm ref}}$ for all $n,j,i$.

The training \emph{loss function} is the $L^1$ error given by,
\begin{equation}
\label{eq:blf}
E_1(w) := \Dx \sum\limits_{i=1}^{I} \sum\limits_{j=1}^J  \sum\limits_{n=1}^{2} |U^{n,i}_j - U^{n,i}_{j,{\rm ref}}|.
\end{equation} 
Here, $w = \{w^{n}\}_{1,2,3}$ for $n=1,2$ represents the vector of $6$ parameters that specifies the scheme \eqref{eq:bfvl} on this grid. We remark that the $L^1$ error is natural for conservation laws as it the norm under which the data to solution operator is continuous \cite{DAF1}.

\begin{table}[htbp]
\centering
\begin{tabular}{|c|c|c|c|c|c|c|c|}
\hline
$w^{1,\ast}_1$ & $w^{1,\ast}_2$ & $w^{1,\ast}_3$ & $w^{2,\ast}_1$  & $w^{2,\ast}_2$ & $w^{2,\ast}_3$ & Gain & Speedup \\
\hline
$0.26$   &  $0.2$   &   $0$ & $0.3$ & $0.3$ & $0$ & $1.42$ & $5.33$ \\ 
\hline
\end{tabular}
\caption{The performance of the trained finite volume Rusanov scheme \eqref{eq:bfvl} on the Burgers' equation \eqref{eq:bg} for the smooth initial data \eqref{eq:kl1}. The optimized weights are labeled by time level (superscript) and spatial location (subscript). Gain is the ratio of the (mean) error with the standard Rusanove scheme and the (mean) error with the trained scheme on the test set. Speedup represents the overall gain in computational efficiency with the trained scheme over the standard Rusanov scheme.}
\protect \label{tab:6}
\end{table}

The loss function is minimized using a simplified version of the stochastic gradient descent algorithm with a batch size of $4$. The stochastic gradient algorithm is applied sequentially i.e, first the minimizers at first time step are determined and then we determine the minimizers at the second time step. This step by step minimization may not yield the optimum in the 6 dimensional parameter space but was found to be reasonably accurate. Moreover, it is computationally cheaper and consistent with the time marching form of numerical methods for evolutionary PDEs. The algorithm was initialized by setting $w^n_{1,2,3} \equiv 0.5$ (corresponding to the standard Rusanov scheme on the coarse grid) and the approximate minimizers, labeled by the vector $w^{\ast}$ are shown in table \ref{tab:6}.  As seen from the table, the minimizing weights are always below the value of $0.5$ (being equal to zero in one case). This implies that the numerical diffusion is reduced during training as the solutions in this case are still identified as mostly smooth. 

\begin{table}[htbp]
\centering
\begin{tabular}{|c|c|c|c|c|c|c|c|}
\hline
$w^{1,\ast}_1$ & $w^{1,\ast}_2$ & $w^{1,\ast}_3$ & $w^{2,\ast}_1$  & $w^{2,\ast}_2$ & $w^{2,\ast}_3$ & Gain & Speedup \\
\hline
$0.24$   &  $0.24$   &   $1.25$ & $0.24$ & $0.26$ & $0$ & $2.48$ & $9.55$ \\ 
\hline
\end{tabular}
\caption{The performance of the trained finite volume Rusanov scheme \eqref{eq:bfvl} on the Burgers' equation \eqref{eq:bg} for the rough initial data \eqref{eq:kl2}. The optimized weights are labeled by time level (superscript) and spatial location (subscript). Gain is the ratio of the (mean) error with the standard Rusanove scheme and the (mean) error with the trained scheme on the test set. Speedup represents the overall gain in computational efficiency with the trained scheme over the standard Rusanov scheme.}
\protect \label{tab:7}
\end{table}

A \emph{test set} is generated by selecting at random $100$ realizations from the initial datum \eqref{eq:kl1} and the \emph{gain}, defined as the ratio of the mean error of the standard Rusanov scheme (on the coarse grid) to the trained scheme is calculated and shown in table \ref{tab:6}. This gain of $1.42$ is rather modest in this case, compared to the previous examples of linear PDEs. However, when we calculate the \emph{speed up} i.e, the ratio of the computational work (time) required to obtain a similar error as the trained scheme (on the coarse grid), but by the standard Rusanov scheme on a finer mesh. In the case of smooth data for the Burgers' equation, this speed up is given in table \ref{tab:6} and amounts to a factor of $5.33$.

\begin{figure}[htbp]
\centering
\includegraphics[width=12cm]{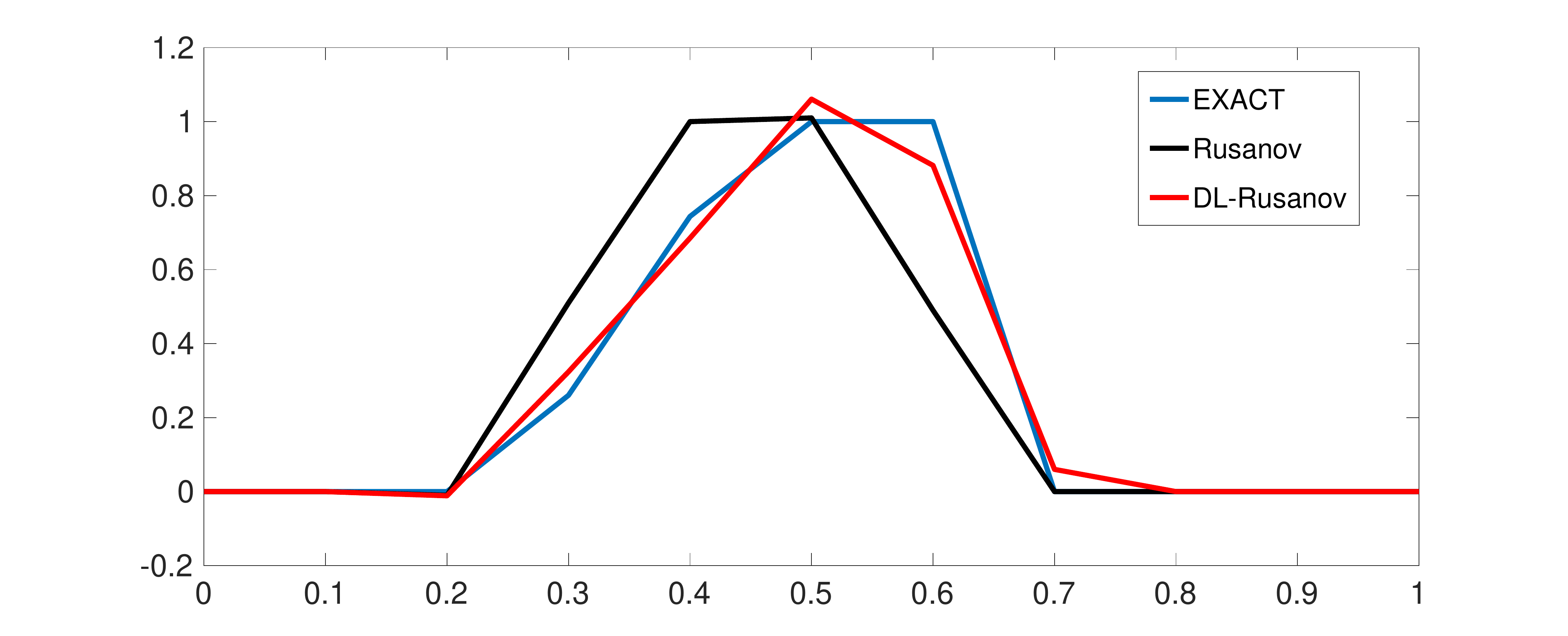}
\caption{Solutions of the Burgers' equation \eqref{eq:bg} at time $T=0.1$ with an initial value from the \emph{test set} for rough data i.e \eqref{eq:kl2} with a particular realization of $Y_{1,2,3}(\omega)$. We show the reference solution (EXACT), solutions approximated with the standard Rusanov scheme and the solutions by computed \emph{Data learned} (DL) trained scheme, \eqref{eq:bfvl} with parameters given in table \ref{tab:6}, on a grid with two time steps and $10$ mesh points.}
\label{fig:8}
\end{figure}

Larger gains are obtained for rough initial data given by the random initial condition \eqref{eq:kl2}. Here, the amplitude of the initial discontinuity and the locations of both jumps are uncertain. In this case, we generate the training data $u^i_0$, identically to the previous case. The reference solution is computed and consists of a right moving shock and a rarefaction on the left (see figure \ref{fig:8}), for each realization. The loss function \eqref{eq:blf} is minimized and the (approximate) minimizers are shown in table \ref{tab:7}. In this case, some optimal values of the weights are significantly higher than $0.5$, indicating larger diffusion around the right moving shock whereas many weights are well below the value of $0.5$, indicating a modulation of numerical diffusion around smooth regions, identified from the training set. The test set is chosen as before and gain, shown in table \ref{tab:7} and amounting to a factor of $2.48$, is higher than the smooth case. Moreover, the overall speed up in this case is $9.55$, representing an order of magnitude computational speed up over the standard Rusanov scheme. 

The solutions computed with the trained scheme and with the standard Rusanov scheme, for one particular realization of the initial data \eqref{eq:kl2} are shown in figure \ref{fig:8}. We observe from the figure that the trained scheme significantly outperforms the Rusanov scheme for this realization and provides an accurate approximation, even on this very coarse grid. 

\section{Euler equations}
\label{sec:7}
The Euler equations of gas dynamics are a prototypical example of a hyperbolic system of conservation laws \cite{DAF1}. In one space dimension, the Euler equations, representing the conservation of mass, momentum and energy are,
\begin{equation}
\label{eq:eul}
\begin{aligned}
u_t + (f(u))_x &= 0, \quad (x,t) \in (0,1) \times (0,T)\\
u &= \left[\rho,\rho v,  E \right], \\
f(u) &= \left[\rho v, \rho v^2 +p, (E+p)v\right], \\
u(x,0) &= u_0(x).
\end{aligned}
\end{equation}
Here $u:(0,1) \times (0,T) \rightarrow \R^3$ is the vector of unknowns and $f:\R^3 \rightarrow \R^3$ is the flux vector. The density of the gas is denoted by $\rho$, the velocity by $v$, the pressure by $p$ and the total energy by $E$. The equations are closed by specified a thermodynamic relation between the variables, such as the \emph{ideal gas equation of state}:
\begin{equation}
\label{eq:eos}
E:= \frac{p}{\gamma - 1} + \frac{1}{2} \rho v^2,
\end{equation}
with gas constant $\gamma$. The system \eqref{eq:eul} is \emph{hyperbolic} with the three eigenvalues,
\begin{equation}
\label{eq:eig}
\lambda_1 = u - a, ~ \lambda_2 = u, \quad \lambda_3 = u+ a,
\end{equation}
given in terms of \emph{sound speed},
\begin{equation}
\label{eq:ss}
a =  \sqrt{\frac{p}{\gamma \rho}}.
\end{equation}

The solutions to systems of conservation laws, such as the Euler equations \eqref{eq:eul}, develop finite time discontinuities such as shock waves and contact discontinuities, even when the initial data is smooth and as in the scalar case, there is a notion of entropy solutions for them.
\subsection{Numerical scheme}
We discrete the computational domain of $[0,1] \times [0,T]$ as in the case of scalar conservation laws (section \ref{sec:clnum}), and evolve cell averages of the vector of unknowns i.e,  $U^n_j = \left[ \rho^n_j, (\rho v)^n_j, E^n_j\right]$ by the (generalized) finite volume scheme,
 \begin{equation}
\label{eq:efv1}
\begin{aligned}
U^{n+1}_j &= U^n_j - \frac{\Dt}{\Dx}\left(F^n_{\jphf} - F^n_{\jmhf} \right), \quad
U_j^n &=  \int_{C_j} u_0(x) dx
\end{aligned}
\end{equation}
Here, $F^n_{\jphf} = F(U^n_j,U^n_{j+1})$ is a numerical flux vector, consistent with the flux vector $f$ in \eqref{eq:eul}. The Local Lax-Friedrichs or \emph{Rusanov} flux for the Euler equations is given by,
\begin{equation}
\label{eq:erus1}
\begin{aligned}
F(U^n_j,U^n_{j+1}) &= \frac{1}{2}\left(f(U^n_j) + f(U^n_{j+1})\right) - \frac{1}{2}\max(|v^n_j|+a^n_j, |v^n_{j+1}| + a^n_{j+1})(U^n_{j+1} - U^n_j),
                                \end{aligned}
\end{equation}
with $a^n_j$ being the sound speed corresponding to the state $U^n_j$, computed from \eqref{eq:ss}. Thus, the numerical diffusion is scaled with an estimate (upper bound) of the local maximal wave speed given in \eqref{eq:eig}.  The primitive variables are calculated from the computed conservative variables at the end of each time step.  The Rusanov flux is known to very diffusive for systems of conservation laws, particularly around contact discontinuities. However, we choose it here in order to be consistent with our choice in the scalar case and to investigate whether we can train a scheme to (significantly) improve on its numerical performance. 

As in the scalar case, we cast the finite volume scheme \eqref{eq:efv1} in our machine learning framework by generalizing the Rusanov flux \eqref{eq:erus1} to
\begin{equation}
\label{eq:erusdl1}
\begin{aligned}
F(U^n_j,U^n_{j+1}) &= \frac{1}{2}\left(f(U^n_j) + f(U^n_{j+1})\right) - \left(w^n_{\jphf} \max(|v^n_j|+a^n_j, |v^n_{j+1}| + a^n_{j+1}) (U^n_{j+1} - U^n_j)\right).
\end{aligned}
\end{equation}
Here $w^n_{\jphf} \in \R,~\forall j$, are weights corresponding to the scaling of the local wave speed. The resulting scheme is given by,
\begin{equation}
\label{eq:efvl}
 \begin{aligned}
U^{n+1}_j &= U^n_j - \frac{\Dt}{2\Dx} \left(f(U^n_{j+1})-f(U^n_{j-1})\right) \\
&+\frac{\Dt}{\Dx} \left(w^n_{\jphf}\max(|v^n_j|+a^n_j, |v^n_{j+1}| + a^n_{j+1})(U^n_{j+1} - U^n_j)\right) \\
&- \frac{\Dt}{\Dx} \left(w^n_{\jmhf}\max(|v^n_j|+a^n_j, |v^n_{j-1}| + a^n_{j-1})(U^n_{j} - U^n_{j-1})\right).
\end{aligned}
\end{equation}
It is straightforward to verify that the scheme \eqref{eq:efvl} is a \emph{conservative} and \emph{consistent} discretization of the one-dimensional Euler equations. It is also (formally) first-order accurate. As in the case of the Burgers' equation, we make a further simplification by \emph{pooling} values of the weights $w^n$ in longer \emph{windows} i.e by requiring \eqref{eq:wind} with $0\leq m \leq J$ as the window length and $j_s \subseteq \{1,\ldots, J\}$ as a subset of grid points on which we center the pooling windows. 

\subsection{Training and results on the test set.}
For our numerical experiments, we consider the one-dimensional Euler equations \eqref{eq:eul} with the ideal gas equation of state \eqref{eq:eos} and gas constant $\gamma = 1.4$, corresponding to a diatomic gas. We fix $\Dx = 0.05$ i.e, we discretize the interval $[0,1]$ into $20$ cells. Our final time is $T=0.15$ and the time interval is divided into \emph{five time steps} with time step size $\Dt = 0.03$. The scheme \eqref{eq:efvl} is closed at the boundary points by imposing \emph{transparent boundary conditions} i.e, a zeroth order extrapolation by  setting $U^n_0 = U^n_1$ and $U^n_{J+1} = U^n_J$.

We also \emph{pool} the weights of the numerical diffusion by a defining a pooling window of size $m=3$ in \eqref{eq:wind}. Hence, at each time step we need to specify $6$ weights namely $\{w^n_{1},w^n_{2},w^n_{3},w^n_{4}, w^n_5,w^n_6\}$ in \eqref{eq:efvl} by grouping every $3$ cell interfaces (starting from the left)., The weights on the boundary interfaces are determined from the transparent boundary conditions. Hence, the scheme \eqref{eq:efvl} contains $30$ parameters that need to be determined in the \emph{training} process.
 
To generate the training set, we consider the following \emph{random} initial data,
\begin{equation}
\label{eq:sst}
\begin{aligned}
\rho_0(x,\omega) &= \begin{cases} 
                           \rho_l + \epsilon Y_1(\omega), &{\rm if} \quad 0 < x < 0.5 + \epsilon Y_2(\omega)  , \\
                            \rho_r + \epsilon Y_3(\omega), &{\rm if} \quad 0.5 + \epsilon Y_2(\omega) < x < 1, 
                           \end{cases} \\
v_0(x,\omega) &= 0, \quad 0 < x < 1, \\
p_0(x,\omega) &= \begin{cases} 
                           p_l + \epsilon Y_4(\omega), &{\rm if} \quad 0 < x < 0.5 + \epsilon Y_2(\omega)  , \\
                            p_r + \epsilon Y_5(\omega), &{\rm if} \quad 0.5 + \epsilon Y_2(\omega) < x < 1, 
                           \end{cases} 
                           \end{aligned}                          
\end{equation}       
Here, $\rho_l = p_l = 1$, $\rho_r = p_r = 0.4$ and $\epsilon = 0.1$. The random variables $Y_{1,2,3,4,5}(\omega)$ are drawn from a uniform distribution on the interval $[-1,1]$. Thus, the initial data \eqref{eq:sst} corresponds to a stochastic version \cite{SSSID1} of the well-known \emph{Sod shock tube} problem for the Euler equations by considering a random interface that separate random jumps in the initial density and pressure. 
 
 We draw $I=50$ samples in \eqref{eq:sst} and label the resulting initial data as $u_0^i$. These data initialize the scheme \eqref{eq:efvl} to generate the updated solutions $U^{n,i}_j$ for all $j$ and $1\leq n \leq 5$. A reference solution is computed using the standard Rusanov scheme (i.e setting $w^n_{\jphf} \equiv 0.5$) in \eqref{eq:efvl} on a fine mesh of $1000$ points and with the time step determined by the corresponding CFL number as in \eqref{eq:cfl1}. This fine grid solution is projected on the underlying coarse grid by averaging. We denote the reference solution as $U^{n,i}_{j,{\rm ref}}$ for all $n,j,i$.
 
The training \emph{loss function} is following version of the $L^1$ error,
\begin{equation}
\label{eq:elf}
E_1(w) := \Dx \sum\limits_{i=1}^{I} \sum\limits_{j=1}^J  \sum\limits_{n=1}^{5} \left(|\rho^{n,i}_j - \rho^{n,i}_{j,{\rm ref}}| + |v^{n,i}_j - v^{n,i}_{j,{\rm ref}}| + |p^{n,i}_j - p^{n,i}_{j,{\rm ref}}| \right)
\end{equation} 
Here, $\rho^{n,i}_j,v^{n,i}_j,p^{n,i}_j$ are the density, velocity and pressure computed on the coarse grid by the numerical scheme \eqref{eq:efvl} for the training initial data $u^{n,i}_0$.  We choose the $L^1$ error in the primitive variables, rather than in the conservative variables. Our choice is motivated by the fact that in practice, one is interested in measuring the velocity and the pressure, rather than the momentum and energy.

We minimize the loss function \eqref{eq:elf} on the $30$-dimensional parameter space by using a stochastic gradient method with batch size of $5$. The stochastic gradient method is initialized by setting all the weights to $w^n_l \equiv 0.5$, corresponding to the standard Rusanov scheme. As in the case of the Burgers' equation, we will optimize the loss function sequentially in time, corresponding to each time step. The stochastic gradient method converges fairly quickly in this case to the optimized weights presented in table \ref{tab:8}. 
\begin{table}[htbp]
\centering
\begin{tabular}{|c|c|c|c|c|c|c|}
\hline
$n$ & $w^{n,\ast}_1$ & $w^{n,\ast}_2$ & $w^{n,\ast}_3$ & $w^{n,\ast}_4$  & $w^{n,\ast}_5$ & $w^{n,\ast}_6$ \\
\hline
$1$   &  $0.5$   &   $0.5$ & $0.42$ & $0.5$ & $0.5$ & $0.5$ \\ 
\hline
$2$   &  $0.5$   &   $0.5$ & $0.23$ & $0.48$ & $0.5$ & $0.5$ \\ 
\hline
$3$   &  $0.5$   &   $-0.3$ & $0.29$ & $0.51$ & $0.5$ & $0.5$ \\ 
\hline
$4$   &  $0.5$   &   $0.17$ & $0.19$ & $0.54$ & $0.5$ & $0.5$ \\ 
\hline
$5$   &  $0.5$   &   $0.2$ & $0.33$ & $0.77$ & $0.3$ & $0.5$ \\ 
\hline
\end{tabular}
\caption{Optimized weights $w^n_l$ for $1\leq n \leq 5$ and $1 \leq l \leq 6$ for the trained Rusanov scheme \eqref{eq:efvl} for the one-dimensional Euler equations \eqref{eq:eul} with the random initial data training set.}
\protect \label{tab:8}
\end{table}

\begin{figure}[htbp]
\centering
\includegraphics[width=0.3\linewidth]{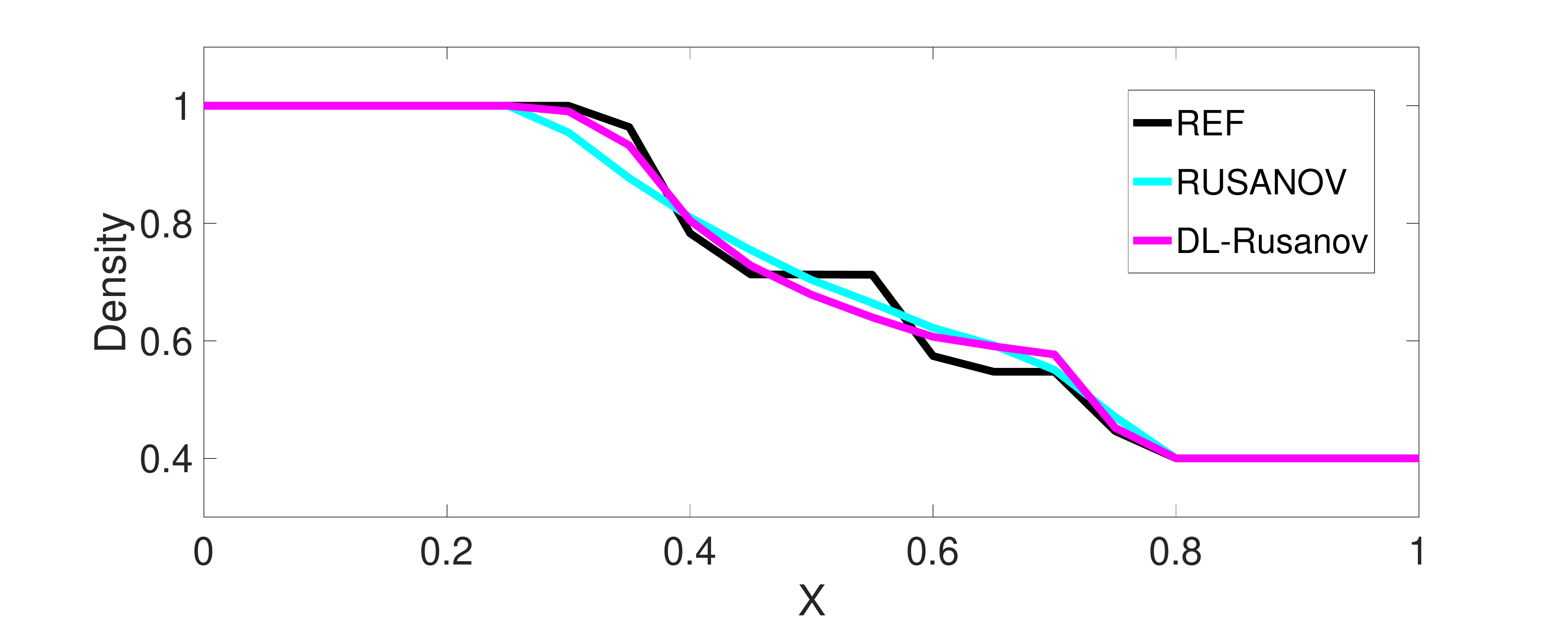} \includegraphics[width=0.3\linewidth]{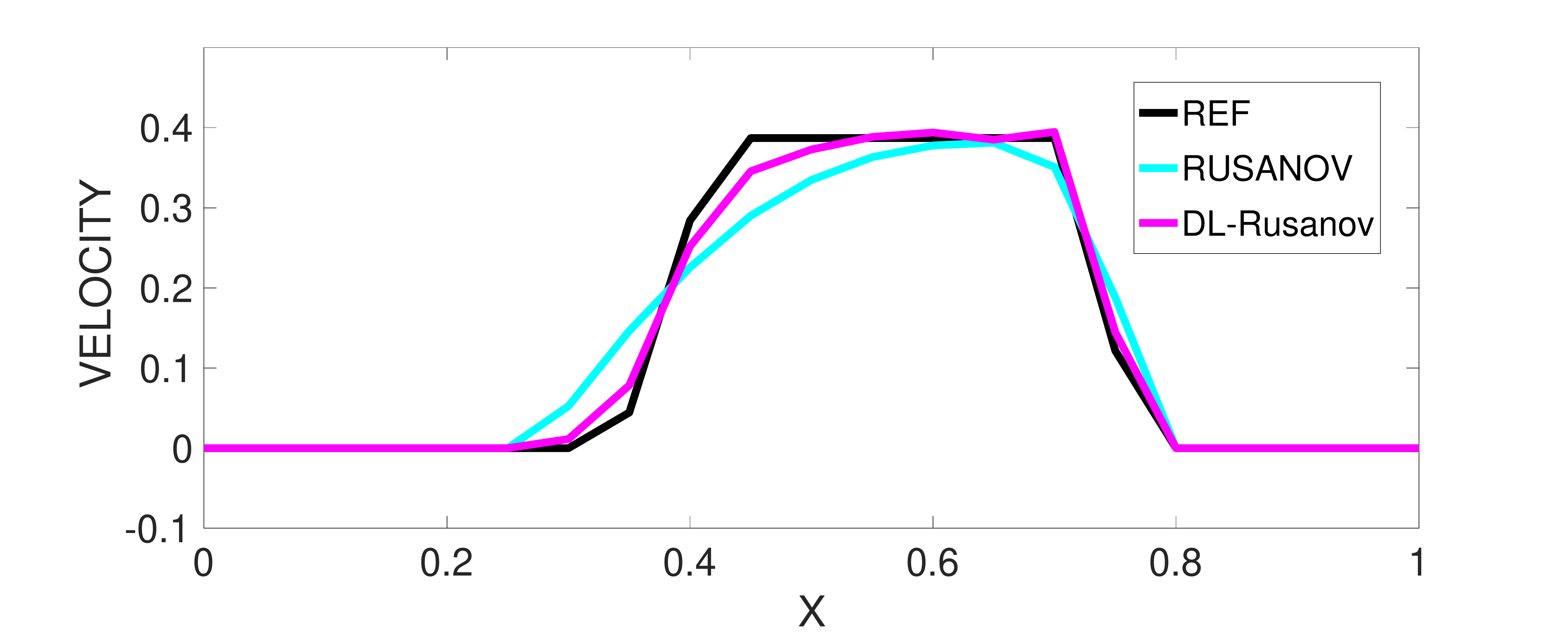}  \includegraphics[width=0.3\linewidth]{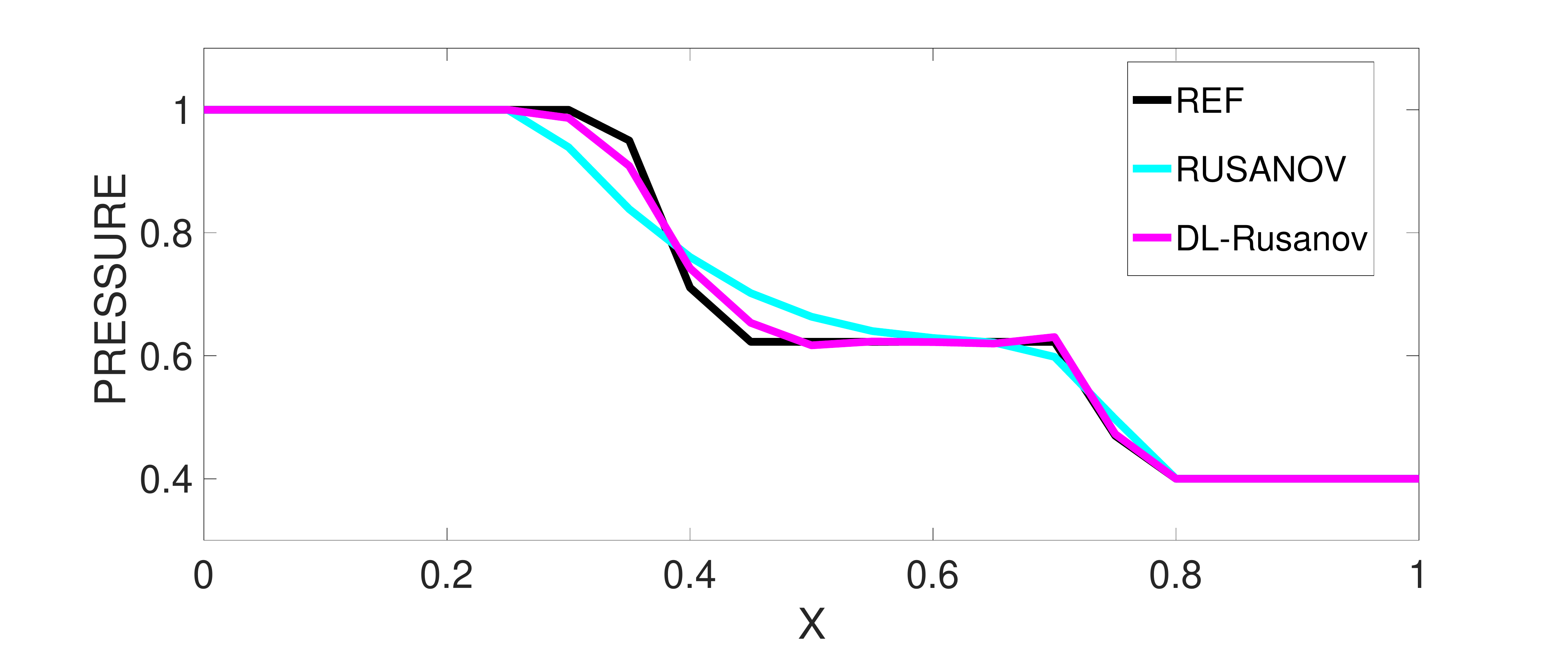} 
\caption{Numerical approximation at time $T=0.15$ of the $1$-dimensional Euler equations \eqref{eq:eul} for a Sod shock tube initial data i.e, \eqref{eq:sst} with $\epsilon = 0$. We compare a reference solution (computed on a fine mesh but initial data on a coarse mesh), the trained scheme \eqref{eq:efvl}, on a grid of $20$ mesh points and $5$ time steps and with optimal weights given in table \ref{tab:8} (DL-Rusanov) and a standard Rusanov scheme on the same coarse grid. Left: Density, Middle: Velocity and Right: Pressure. Note the much improved approximation of the shock wave and the rarefaction wave and further diffusion of the contact discontinuity with the trained scheme.}
\label{fig:9}
\end{figure}

The (approximate) optimized weights, shown in table \ref{tab:8}, follow an interesting pattern. We recall that the solutions of the Euler equations \eqref{eq:eul} with initial data \eqref{eq:sst} consist of the initial discontinuity breaking down
into three waves namely, a left moving rarefaction, a right moving contact and an even faster right moving shock wave, see figure \ref{fig:9} for a snapshot of the solution. First, we observe that only $13$ of the $30$ weights assume a value different from the initial value and this variation is more pronounced with time as waves develop, separate and move away from the initial discontinuity. The values assumed by the weights include locations where the optimal weights are greater than $0.5$, implies that more diffusion is added but there are many locations, particularly on the left of the initial interface, corresponding to the rarefaction wave, where the weights are significantly less than the initial value of $0.5$ (one is even negative), indicating that diffusion is removed near the continuous part of the solution.

We generate a \emph{test set} by drawing $1000$ samples from \eqref{eq:sst}. The resulting gain, defined as the ratio of the (mean) error of the standard Rusanov scheme on the test set, to the (mean) error of the trained scheme, is computed and is determined to be a factor of $2.17$. Although this seems rather modest given the much more significant gains for the heat and the linear transport equation, it is comparable to the gain for the Burgers' equation reported in table \ref{tab:7}. However, the key quantity to demonstrate the efficiency of the trained scheme is the \emph{speed up} i.e the ratio of the (mean) computational time for the standard Rusanov scheme (on a finer grid) to achieve the same error as the trained scheme on the underlying coarse grid. Given that the observed (mean) order of convergence of the standard Rusanov scheme on this problem is found to be $0.57$ (even if the Rusanov scheme is (formally) first-order accurate), we obtain that the trained scheme on a grid of $20$ mesh points (and five time steps) is comparable in error to a standard Rusanov scheme on a grid of $80$ mesh points (and time steps determined from the CFL number). Consequently, the speedup is a factor of $16$.

Further insight into the performance of the trained scheme is provided in figure \ref{fig:9}, where we plot the computed density, velocity and pressure with the trained scheme \eqref{eq:efvl} and compare it with the standard Rusanov scheme and a reference solution. As seen from the figure \ref{fig:9} (left), the trained scheme provides a considerably more accurate approximation of the rarefaction wave and the (fast) shock, even on this very coarse grid. On the other hand, it dissipates the contact even further. This counter-intuitive behavior can be explained on the basis of the loss function \eqref{eq:elf}. Note that the pressure and the velocity are constant across the contact. Thus, the contribution of the contact in the loss function \eqref{eq:elf} can be rather small. Hence, during the training process, the scheme "decides" not to approximate the contact better but to focus on approximating the shock and the rarefaction more accurately. This strategy clearly bears fruit as the velocity (figure \ref{fig:9} (middle)) and pressure (figure \ref{fig:9} (right)) are approximated very well, albeit with small oscillations, leading to a larger overall reduction in the loss function \eqref{eq:elf}. This non-intuitive behavior is in stark contrast to traditional approaches that focus on approximating the contact discontinuity more accurately.

\section{Discussion}
\label{sec:8}
Numerical methods for efficient approximation of (time-dependent) ordinary and partial differential equations are well established. However, emerging applications such as uncertainty quantification (UQ), (Bayesian) inverse problems and
(real time) optimal control and design require fast (computationally cheap) yet accurate numerical methods. Existing numerical schemes, particularly for nonlinear PDEs, fail to provide reasonable accuracy at very low computational cost.

In this paper, we have proposed a \emph{machine learning} framework, summarized in Algorithm \ref{alg:1}, for designing such cheap yet accurate methods. The basis of our algorithm is the observation that the computational cost of existing numerical methods on (very) coarse (space-time) grids is rather low. However, these methods are too inaccurate on such grids to be of practical use. We aim to increase the accuracy (reduce the numerical error) of these methods on coarse grids. To this end, we recast of generalizations of standard numerical methods in terms of artificial neural networks i.e layers of units coupled with linear operators and possibly nonlinear activations, but with a set of undetermined parameters. These parameters are \emph{trained} to minimize a \emph{loss function} on a carefully chosen \emph{training set} in an offline training phase. 
The key properties of our proposed algorithm are
\begin{itemize}
\item The resulting method is always consistent with underlying ODE or PDE by design. Additional constraints can be imposed on the trainable parameters to ensure stability.
\item The method is guaranteed to be more accurate than a standard numerical method on the same grid as the gradient descent method is initialized with parameters that correspond to a standard method.
\item The method is very simple to implement with minor changes in existing numerical ODE and PDE solvers.
\end{itemize}
Although no theoretical guarantees have been established on whether the proposed algorithm significantly outperforms standard methods on the \emph{test set}, we have presented extensive numerical experiments to ascertain this enhancement in performance. Our numerical experiments include a linear and a non-linear ODE, the linear heat and transport equations, scalar conservation laws (Burgers' equation) and the Euler equations of gas dynamics. We considered underlying numerical schemes that include implicit multi-step methods for ODEs, implicit finite difference schemes for linear PDEs and explicit finite volume schemes for nonlinear PDEs. Loss functions, measuring error in either $L^1$ or $L^2$ norms were minimized with stochastic gradient methods. The numerical experiments demonstrated a significant gain in performance (computational speed up) over the underlying standard numerical method. The gains ranged from an order of magnitude for nonlinear problems to two (or three) orders of magnitude for linear problems. In all cases considered here, the machine learning algorithm \ref{alg:1} provided a numerical method with reasonable accuracy on a very coarse space-time grid. Hence, it could serve as a basis for the solution of complex problems in UQ, inverse problems or (real time) optimal control. 

It is instructive to compare our approach to possible \emph{deep learning} of the solutions of differential equations. It is essential to recall that one can cast standard numerical methods for time-dependent differential equations in the form of a deep neural network, see section \ref{sec:2} (figure \ref{fig:2}) and section \ref{sec:6} (figure \ref{fig:7}) for concrete examples. At the very least, non-linearities that occur in standard numerical methods, can be approximated with standard deep neural networks based on ReLU activations. Thus, our approach does consist of approximating differential equations with a form of deep networks. Hence, standard deep learning methodology such as back propagation, stochastic gradients, pooling etc and software frameworks like \emph{TENSORFLOW}, can be readily used.

However, as explained in section \ref{sec:rusnn}, there is a significant difference in our algorithm with customary deep learning. In machine learning, one usually trains all the parameters in the network. Doing so in our context, see figure \ref{fig:7}, may lead to a lack of consistency with the underlying differential equation. In order to retain consistency (and possibly stability), one needs to constrain the set of parameters in order to recover these properties for every value of the trainable parameters. We do so with our (natural) generalization of numerical schemes. Thus, one can consider algorithm \ref{alg:1} as a deep learning algorithm, with a very particular architecture, and with a (very) restricted set of trainable parameters. We retain consistency and at the same time, notice significant gains in computational efficiency. It could be that a free training of all the parameters in the deep network, underlying our generalized numerical method, will automatically identify regions of the parameter space (particularly if additional penalization terms are added to the loss function) to retain consistency and stability. This approach needs to be explored in the future. 

It should be emphasized that the (approximate) optimal values of parameters calculated in all our examples, depend strongly on the training set, the underlying coarse grid and parameters of the problem such as the diffusion coefficient for the heat equation \eqref{eq:ht} or the wave speed in the linear advection equation \eqref{eq:lt}. One can also train the algorithm to take a possible stochastic model for some of these parameters into account. Moreover for sake of simplicity of exposition, we have mostly considered model problems with very simple underlying numerical methods in this paper. The number of undetermined parameters was in the range of $2-3$ for ODEs and linear PDEs to at most $30$ for nonlinear PDEs. State of the art deep learning architectures handle hundreds of thousands to millions of trainable parameters and we anticipate a much larger gain in efficiency when we dramatically increase the width and depth of our schemes, represented as networks. Applications of the proposed algorithm \ref{alg:1} to realistic multi-dimensional problems in uncertainty quantification and inverse problems is the subject of ongoing work.   
\section*{Acknowledgements}
The author thanks Kjetil O. Lye (SAM, ETH Z\"urich) for interesting discussions on this topic. The research of SM was partially funded by ERC CoG 770880 COMANFLO.

\end{document}